\documentclass[11pt]{article}
\usepackage{amssymb,amsmath,fullpage}
\title{\bf Pad\'e interpolation and hypergeometric series}
\author{Masatoshi Noumi \\
{\normalsize Department of Mathematics, Kobe University}
}
\date{} 
\newtheorem{thm}{Theorem}[section]

\newtheorem{lem}[thm]{Lemma}

\numberwithin{equation}{section}
\newcommand{\dsum}[2]{{\displaystyle\sum_{#1}^{#2}}}

\newcommand{\dprod}[2]{{\displaystyle\prod_{#1}^{#2}}}
\newcommand{\tprod}[2]{{\textstyle\prod_{#1}^{#2}}}
\newcommand{\proof}[2]{
\par\noindent{\sl Proof {#1}$:$}\ 
{#2}\hfill$\square$\par\medskip
}
\begin{document}
\maketitle
\begin{quote}
{\bf Abstract: }  We propose a class of Pad\'e interpolation problems whose 
general solution is expressible in terms of determinants of hypergeometric series. 
\end{quote}
%%%%%%%
\baselineskip=16pt
%%%%%%%
\section{Introduction}

In this paper we investigate a class of Pad\'e interpolation 
problems to which the solutions are expressible in terms of determinants 
of hypergeometric series. 
Pad\'e interpolation problems have been discussed 
by Spiridonov--Zhedanov \cite{SZ2007}
from the viewpoint of biorthogonal rational functions. 
They are also sources of the Lax pairs for discrete Painlev\'e equations constructed 
by Yamada \cite{Y2009FE}, \cite{Y2009Sigma}, and by Noumi--Tsujimoto--Yamada 
\cite{NTY2013}.  
The goal of this paper is to clarify how hypergeometric series arise % widely 
in Pad\'e interpolation problems, by analyzing the determinantal expression 
of the general solution. 

In Section 2, we formulate a general Pad\'e interpolation problem and 
a universal determinant formula for the general solution 
(Theorem \ref{thm:PQGen}).  
We also show that the determinants expressing the general solution 
can be condensed to smaller determinants 
by a variation of the Dodgson condensation
(Theorem \ref{thm:PQUV}).  
After these preliminaries, we investigate in Section 3 
a class of Pad\'e interpolation problems relevant to 
generalized hypergeometric series 
%%% added in revision
${}_{r+1}F_{r}$.  
%%%%%%%%%%%%%%
We propose there two types of formulas expressing the 
solutions in terms of determinants of generalized hypergeometric series.  
The first one (Theorem \ref{thm:HG}), 
derived through Theorem \ref{thm:PQUV}, is based on 
the condensation of determinants 
and Krattenthaler's determinant formula, 
while the second (Theorem \ref{thm:HG2}) is constructed 
by means of the Saalsch\"utz summation 
formula for terminating ${}_{3}F_2$ series.  
%%% added in revision
We remark that Pad\'e approximations to 
generalized hypergeometric functions  
have been discussed by Luke \cite{L1975}, \cite{L1977}.  
It would be an important question to clarify the relationship 
between interpolations and approximations in the context of 
generalized hypergeometric functions. 
%%%%%%%%%%%%%%
Section 4 is devoted to the extension of these results to 
three types of 
{very well-poised hypergeometric series} including basic (trigonometric) and elliptic 
hypergeometric series.  
The two determinant formulas of Theorem \ref{thm:EHG} and Theorem \ref{thm:EHG2}
are obtained by Warnaar's elliptic extension of the Krattenthaler determinant 
and by the Frenkel--Turaev summation formula for terminating ${}_{10}V_{9}$ series. 

Two fundamental tools of our approach are the condensation of 
determinants {\em along a moving core} 
(an identity of Sylvester type), and  
variations of Krattenthaler's determinant formula. 
For the sake of convenience, 
these subjects are discussed separately in Appendix A and 
Appendix B respectively.  
%%%%% added in revision
Generalization of Sylvester's identity on determinants 
has been developed extensively by M\"uhlbach--Gasca \cite{MG1985}
(see also \cite{AAM1996}).  
The version we use in this paper (Lemma \ref{lem:ACond2}), 
based on the {\em Neville elimination strategy},
is originally due to Gasca--L\'opez-Carmona--Ramirez \cite{GLR1982}. 
We also remark that Sylvester's identity and its extensions 
play important roles in recent studies of integrable systems
(see Spicer--Nijhoff--van der Kamp \cite{SNK2011} for example). 
As to Appendix B, basic references are the works of Krattenthaler 
\cite{K1999}, \cite{K2005} and Warnaar \cite{W2002}
(see also Normand \cite{N2004} for recent works on the evaluation 
of determinants involving shifted factorials).  
Although the contents of these appendices are basically found  
in the literature, we include them as self-contained expositions 
which might be helpful to the reader.  

\par\medskip
Throughout this paper we use the following notation of submatrices and minor 
determinants.  
For an $m\times n$ matrix 
$X=\big(x_{ij}\big)_{1\le i\le m, 1\le j\le n}$ (with entries in a commutative ring), 
we denote by 
\begin{align}
\arraycolsep=2pt
X^{i_1,\ldots,i_r}_{j_1,\ldots,j_s}
=\left[\,
\begin{matrix}
x_{i_1j_1} & \ldots & x_{i_1j_s}\\[4pt]
\vdots &\ddots& \vdots \\
x_{i_rj_1} & \ldots & x_{i_rj_s}
\end{matrix}\,
\right]
=
\big(x_{i_aj_b}\big)_{1\le a\le r,1\le b\le s}
\end{align}
the $r\times s$ submatrix of $X$ with row 
indices $i_1,\ldots,i_r\in\{1,\ldots,m\}$ and 
column indices $j_1,\ldots,j_r\in\{1,\ldots,n\}$.  
When $r=s$, we denote by 
$\det X^{i_1,\ldots,i_r}_{j_1,\ldots,j_r}$  
the corresponding minor determinant.

\section{Pad\'e interpolation problems and their determinant solutions} 
In this section 
%We first 
we formulate general Pad\'e interpolation problems and propose 
some universal determinant formulas for the solutions. 
\par\medskip
Let 
$f_0(x)$, $f_1(x)$, $\ldots$, $f_m(x)$ and 
$g_0(x)$, $g_1(x)$, $\ldots$, $g_n(x)$ be two sequences 
of linearly independent meromorphic functions in 
$x\in\mathbb{C}$ and set $N=m+n$. 
We consider a pair $(P_m(x),Q_n(x))$ of two functions
\begin{align}
P_m(x)&=p_{m0}\,f_0(x)+p_{m1}\,f_1(x)+\cdots+p_{mm}\,f_m(x),
\nonumber\\[4pt]
Q_n(x)&=\,q_{n0}\,g_0(x)\,+\,q_{n1}\,g_1(x)\,+\,\cdots\,+\,q_{nn}\,g_n(x), 
\end{align}
which are expressed as $\mathbb{C}$-linear combinations of $f_j(x)$ and $g_j(x)$ respectively. 
Noting that the ratio $P_m(x)/Q_n(x)$ contains $N+1=m+n+1$ arbitrary constants, 
we investigate the interpolation problem
\begin{align}
\dfrac{P_m(u_0)}{Q_n(u_0)}=v_0,\ 
\dfrac{P_m(u_1)}{Q_n(u_1)}=v_1,\ 
\ldots,\ 
\dfrac{P_m(u_N)}{Q_n(u_N)}=v_N
\end{align}
for a set of $N+1$ generic reference points $x=u_0,u_1,\ldots,u_N$ and 
a set of $N+1$ prescribed values $v_0,v_1,\ldots,v_N$.  
This problem is equivalently rewritten as  
\begin{align}\label{eq:PIlambdamu}
P_m(u_k): Q_n(u_k)=\lambda_k: \mu_k\qquad(k=0,1,\ldots,N)
\end{align}
for $v_{k}=\lambda_k/\mu_k$ ($k=0,1,\ldots,N$). 
We remark that the {\em Pad\'e interpolation problem} 
defined as above contains the {\em Lagrange interpolation problem} as 
a special case where $n=0$ and $g_0(x)=1$.  

\par\medskip
A general solution of this Pad\'e interpolation problem is  
given as follows in terms of $(N+2)\times(N+2)$ determinants: 

\begin{align}
\label{eq:PDet}
\arraycolsep=2pt
P_m(x)
&=\det
\left[\,
\begin{matrix}
f_0(x) & \cdots & f_m(x) &
0 & \cdots & 0 \\[4pt]
\mu_0 f_0(u_0) & \cdots & \mu_0 f_m(u_0) &
\lambda_0 g_0(u_0) & \cdots & \lambda_0 g_n(u_0) \\
\vdots &&\vdots&\vdots&&\vdots\\
\mu_N f_0(u_N) & \cdots & \mu_N f_m(u_N) &
\lambda_N g_0(u_N) & \cdots & \lambda_N g_n(u_N)
\end{matrix}\,
\right],
\\[10pt]
\label{eq:QDet}
Q_n(x)
&=
-\det
\left[\,
\begin{matrix}
0 & \cdots & 0 &
g_0(x) & \cdots & g_n(x) \\[4pt]
\mu_0 f_0(u_0) & \cdots & \mu_0 f_m(u_0) &
\lambda_0 g_0(u_0) & \cdots & \lambda_0 g_n(u_0) \\
\vdots &&\vdots&\vdots&&\vdots\\
\mu_N f_0(u_N) & \cdots & \mu_N f_m(u_N) &
\lambda_N g_0(u_N) & \cdots & \lambda_N g_n(u_N)
\end{matrix}\,
\right].  
\end{align}
\par\medskip
\begin{thm}\label{thm:PQGen}
The pair of functions $(P_m(x),Q_n(x))$ defined by \eqref{eq:PDet}, 
\eqref{eq:QDet}
solves the Pad\'e interpolation problem \eqref{eq:PIlambdamu} 
if $(P_m(u_k),Q_n(u_k))\ne(0,0)$ for $k=0,1,\ldots,N$.  
\end{thm}

In order to prove that this pair $(P_m(x),Q_n(x))$ gives a solution of 
the interpolation problem, we introduce two parameters 
$\lambda,\mu$ and consider 
the $(N+2)\times(N+2)$ determinant 
\begin{equation}
\arraycolsep=2pt
R_{m,n}(x;\lambda,\mu)
=\det
\left[
\begin{matrix}
\mu\, f_0(x) & \cdots & \mu\,f_m(x) &
\lambda\,g_0(x) & \cdots & \lambda\,g_n(x) \\[4pt]
\mu_0 f_0(u_0) & \cdots & \mu_0 f_m(u_0) &
\lambda_0 g_0(u_0) & \cdots & \lambda_0 g_n(u_0) \\
\vdots &&\vdots&\vdots&&\vdots\\
\,\mu_N f_0(u_N) & \cdots & \mu_N f_m(u_N) &
\lambda_N g_0(u_N) & \cdots & \lambda_N g_n(u_N)\, 
\end{matrix}
\right]. 
\end{equation}
By decomposing the top row as 
\begin{equation}
\mu\,(f_0(x),\cdots,f_m(x),\ 0,\cdots,0)
+
\lambda
(0,\cdots,0,\ 
g_0(x),\cdots,g_n(x) ), 
\end{equation}
we have 
\begin{equation}
R_{m,n}(x;\lambda,\mu)=\mu\,P_m(x)-\lambda\,Q_n(x). 
\end{equation}
On the other hand, 
the determinantal expression of $R_{m,n}(x;\lambda,\mu)$ implies 
\begin{equation}
R_{m,n}(u_k;\lambda_k,\mu_k)=\mu_k P_m(u_k)-\lambda_k Q_n(u_k)=0
\qquad(k=0,1,2\ldots,N), 
\end{equation}
and hence
\begin{equation}
P_m(u_k):Q_n(u_k)  = \lambda_k :\mu_k\qquad(k=0,1,\ldots,N)
\end{equation}
as desired.

\par\medskip
The $(N+2)\times(N+2)$ determinants \eqref{eq:PDet}, \eqref{eq:QDet}
representing $P_m(x)$ and $Q_n(x)$ can be 
{\em condensed\/} into an $(m+1)\times(m+1)$ and 
$(n+1)\times(n+1)$ determinants respectively, by means 
of a variation of the Dodgson condensation 
(see Appendix A).  

We denote by 
\begin{align}
F=\big(f_j(u_i)\big)_{0\le i\le N,\ 0\le j\le m},\qquad
G=\big(g_j(u_i)\big)_{0\le i\le N,\ 0\le j\le n}
\end{align}
the matrices defined by the values of the functions 
$f_j(x)$ ($0\le j\le m$) and $g_j(x)$ ($0\le j\le n$), respectively, 
at the reference points $u_i$ ($0\le i\le N$). 
We assume that the configuration of 
reference points $u_k$ $(k=0,1,\ldots,N)$ is {\em generic} in the sense that 
the minor determinants
\begin{align}
\det F^{i,\ldots,i+n}_{0,\ldots,n}\quad(0\le i\le m),\quad
\det G^{i,\ldots,i+m}_{0,\ldots,m}\quad(0\le i\le n)
\end{align}
of maximal size with consecutive rows are all nonzero.  

Assuming that $\lambda_k\ne 0$, $\mu_k\ne0$ for 
$k=0,1,\ldots,N$, 
we set 
\begin{equation}\label{eq:defU}
\arraycolsep=2pt
U_{i,j}
=
\dfrac{\lambda_i}{\mu_i}
\det\left[\,
\begin{matrix}
\tfrac{\mu_i}{\lambda_i}\,f_j(u_i) & g_0(u_i) & \ldots & g_n(u_i)\\[4pt]
\tfrac{\mu_{i+1}}{\lambda_{i+1}}\,f_j(u_{i+1}) & g_0(u_{i+1}) 
& \ldots & g_n(u_{i+1})\\
\vdots & \vdots & & \vdots \\
\tfrac{\mu_{i+n+1}}{\lambda_{i+n+1}} f_j(u_{i+n+1}) & g_0(u_{i+n+1}) & \ldots 
& g_n(u_{i+n+1})
\end{matrix}
\,\right]\,
\big(\det G^{i+1,\ldots,i+n+1}_{0,1,\ldots,n}\big)^{-1}
\end{equation}
for $0\le i<m$, $0\le j\le m $ and
\begin{equation}\label{eq:defV}
\arraycolsep=2pt
V_{i,j}
=
\dfrac{\mu_i}{\lambda_i}\,
\det\left[\,
\begin{matrix}
\tfrac{\lambda_i}{\mu_i}\,g_j(u_i) & f_0(u_i) & \ldots & f_m(u_i)\\[4pt]
\tfrac{\lambda_{i+1}}{\mu_{i+1}}\,g_j(u_{i+1}) & f_0(u_{i+1}) & \ldots & f_m(u_{i+1})\\
\vdots & \vdots & & \vdots \\
\tfrac{\lambda_{i+m+1}}{\mu_{i+m+1}} g_j(u_{i+m+1}) & f_0(u_{i+m+1}) & \ldots 
& f_m(u_{i+m+1})
\end{matrix}
\,\right]\,
\big(\det F^{i+1,\ldots,i+m+1}_{0,1,\ldots,m}\big)^{-1}
\end{equation}
for $0\le i<n$, $0\le j\le n$.
Then by Lemma \ref{lem:ACond2} (of {\em condensation along a moving core}), 
the $(N+2)\times(N+2)$ determinants \eqref{eq:PDet}, \eqref{eq:QDet} 
are condensed as follows 
into $(m+1)\times(m+1)$ and $(n+1)\times(n+1)$ determinants 
respectively (see also \eqref{eq:A18}, \eqref{eq:A19}). 

\begin{thm}\label{thm:PQUV}
The two functions $P_m(x)$, $Q_n(x)$ defined in Theorem \ref{thm:PQGen} 
are expressed as follows in terms of $(m+1)\times(m+1)$ and 
$(n+1)\times(n+1)$ determinants respectively\,$:$ 
\begin{align}\label{eq:PDetU}
\arraycolsep=2pt
P_m(x)&=
\dprod{i=0}{m-1}\,\mu_i\,\dprod{i=0}{n}\,\lambda_{m+i}\,
\det G{\,}^{m,\ldots,m+n}_{0,1,\ldots,n}
\det\left[\,
\begin{matrix}
f_0(x) & \ldots & f_m(x) \\[4pt]
U_{0,0} & \ldots & U_{0,m} \\
\vdots &&\vdots\\
U_{m-1,0}&\ldots & U_{m-1,m}
\end{matrix}\,
\right], 
\\
\label{eq:QDetV}
Q_n(x)&=
\epsilon_{m,n}\dprod{i=0}{n-1}\lambda_i\dprod{i=0}{m}\mu_{n+i}
\det F{\,}^{n,\ldots,n+m}_{0,1,\ldots,m}
\det\left[\,
\begin{matrix}
g_0(x) & \ldots & g_n(x) \\[4pt]
V_{0,0} & \ldots & V_{0,n} \\
\vdots &&\vdots\\
V_{n-1,0}&\ldots & U_{n-1,n}
\end{matrix}\,
\right], 
\end{align}
where $\epsilon_{m,n}=(-1)^{mn+m+n}$. 
\end{thm}

By expanding the determinants $U_{ij}$ and $V_{ij}$ 
along the first column we further obtain the series expansions
\begin{align}\label{eq:SeriesU}
U_{ij}&=\dsum{k=0}{n+1}(-1)^k 
\dfrac{\mu_{i+k}\,\lambda_i}{\lambda_{i+k}\,\mu_i} 
f_j(u_{i+k})
\,\dfrac{\det G^{i,\ldots,\widehat{i+k},\ldots,i+n+1}_{0,1,\ldots,n}}
{\det G^{i+1,\ldots,i+n+1}_{0,1,\ldots,n}}
\qquad(0\le i<m, \ 0\le j\le m),
\\
\label{eq:SeriesV}
V_{ij}&=\dsum{k=0}{m+1}(-1)^k 
\dfrac{\lambda_{i+k}\,\mu_i}{\mu_{i+k}\,\lambda_i} 
g_j(u_{i+k})
\,\dfrac
{\det F^{i,\ldots,\widehat{i+k},\ldots,i+m+1}_{0,1,\ldots,m}}
{\det F^{i+1,\ldots,i+m+1}_{0,1,\ldots,m}}
\qquad(0\le i<n, \ 0\le j\le n).
\end{align}
Hence the problem to determine $P_m(x)$ and $Q_n(x)$ 
is reduced to the computation of minor determinants of the 
matrices $F=\left(f_j(u_i)\right)_{i,j}$ and  $G=\left(g_j(u_i)\right)_{i,j}$. 
We remark that these formulas for $P_m(x)$ and $Q_n(x)$ 
hold {\em universally} 
for any choice of the functions $f_j(x)$ and $g_j(x)$.  

\par\medskip
In Sections 3 and 4, we show that these expansion formulas 
\eqref{eq:SeriesU}, \eqref{eq:SeriesV} in fact 
give rise to {\em hypergeometric series} of various types
for appropriate choices of the functions $f_j(x)$, $g_j(x)$, 
the reference points $u_k$ and the prescribed values 
$v_k=\lambda_k/\mu_k$. 

\section{Hypergeometric series arising from determinants}

We explain below how the series expansions 
\eqref{eq:SeriesU}, \eqref{eq:SeriesV} 
can be used for generating hypergeometric series. 
In this section we use the notation of shifted factorials
\begin{align}
(a)_n=a(a+1)\cdots(a+n-1)=\Gamma(a+n)/\Gamma(a)\qquad (n=0,1,2,\ldots).
\end{align}
\par\medskip
As a typical example, 
we consider the two sequences of rational functions
\begin{equation}\label{eq:funcfg}
f_j(x)=\dfrac{(a+x)_j}{(b+x)_j},\qquad 
g_j(x)=\dfrac{(c+x)_j}{(d+x)_j}\qquad(j=0,1,2,\ldots)
\end{equation}
with four complex parameters $a$, $b$, $c$, $d$, to form  
a pair $(P_m(x),Q_n(x))$ of rational functions 
\begin{equation}\label{eq:PadeR0}
P_m(x)=\dsum{j=0}{m}\ p_{m,j}\dfrac{(a+x)_j}{(b+x)_j},
\qquad
Q_n(x)=\dsum{j=0}{n}\ q_{n,j}\dfrac{(c+x)_j}{(d+x)_j}. 
\end{equation}
Taking an arithmetic progression $u_k=u+k$ 
($k=0,1,2,\ldots,N$; $N=m+n$) of points in $\mathbb{C}$,  
we consider the Pad\'e interpolation problem 
\begin{align}\label{eq:PadeR1}
\dfrac{P_m(u+k)}{Q_n(u+k)}=\dfrac{\lambda_k}{\mu_k}\qquad(k=0,1,\ldots,N).
\end{align}

\begin{thm}\label{thm:HG}
Consider the Pad\'e interpolation problem \eqref{eq:PadeR0}, \eqref{eq:PadeR1}
for the rational functions 
$f_j(x)$, $g_j(z)$ in \eqref{eq:funcfg} 
and the reference points $u_k=u+k$ $(k=0,1,\ldots,N; N=m+n)$. 
Then the solution $(P_m(x),Q_n(x))$ of Theorem 
\ref{thm:PQGen} is explicitly given by 
\begin{align}\label{eq:PDetUHG}
\arraycolsep=2pt
P_m(x)&=
\dfrac{
\tprod{i=1}{n}\,i!\,(d-c)_i
}{
\tprod{i=0}{n}\,(d+u_{m+i})_n
}
\dprod{i=0}{m-1}\,\mu_i\,\dprod{i=0}{n}\,\lambda_{m+i}\,
\det\left[\,
\begin{matrix}
f_0(x) & \ldots & f_m(x) \\[4pt]
U_{0,0} & \ldots & U_{0,m} \\
\vdots &&\vdots\\
U_{m-1,0}&\ldots & U_{m-1,m}
\end{matrix}\,
\right], 
\\
\label{eq:QDetVHG}
Q_n(x)&=
\epsilon_{m,n}
\dfrac{
\tprod{i=1}{m}\,i!\,(b-a)_i
}{
\tprod{i=0}{m}\,(b+u_{n+i})_m
}
\dprod{i=0}{n-1}\lambda_i\dprod{i=0}{m}\mu_{n+i}
\det\left[\,
\begin{matrix}
g_0(x) & \ldots & g_n(x) \\[4pt]
V_{0,0} & \ldots & V_{0,n} \\
\vdots &&\vdots\\
V_{n-1,0}&\ldots & U_{n-1,n}
\end{matrix}\,
\right], 
\end{align}
where 
\begin{align}
U_{ij}&=
\dfrac{(a+u_i)_j}{(b+u_i)_j}
\dsum{k=0}{n+1}
\dfrac{(-n-1)_k(d+u_{i+n})_k}{k!\ (d+u_i)_k}
\dfrac{(a+u_{i+j})_k(b+u_i)_k}{(a+u_{i})_k(b+u_{i+j})_k}
\dfrac{\mu_{i+k}\,\lambda_{i}}{\lambda_{i+k}\,\mu_i}, 
\\
V_{ij}&=
\dfrac{(c+u_i)_j}{(d+u_i)_j}
\dsum{k=0}{m+1}
\dfrac{(-m-1)_k(b+u_{i+m})_k}{k!\ (b+u_i)_k}
\dfrac{(c+u_{i+j})_k(d+u_i)_k}{(c+u_{i})_k(d+u_{i+j})_k}
\dfrac{\lambda_{i+k}\,\mu_i}{\mu_{i+k}\,\lambda_{i}}. 
\end{align}  
\end{thm}

As we remarked in the previous section, 
the functions $U_{ij}$ in Theorem \ref{thm:PQUV} are expressed as 
\begin{align}
U_{ij}&=\dsum{k=0}{n+1}(-1)^k 
\dfrac{\mu_{i+k}\,\lambda_i}{\lambda_{i+k}\,\mu_i} 
f_j(u_{i+k})
\,\dfrac{\det G^{i,\ldots,\widehat{i+k},\ldots,i+n+1}_{0,1,\ldots,n}}
{\det G^{i+1,\ldots,i+n+1}_{0,1,\ldots,n}}. 
\end{align}
Since $u_k=u+k$ $(k=0,1,\ldots,N)$, we have 
\begin{align}
f_j(u_{i+k})=
\dfrac{(a+u_{i+k})_j}{(b+u_{i+k})_j}
%=\dfrac{(a+u_{i})_j}{(b+u_{i})_j}
%\dfrac{(a+u_{i+k})_j(b+u_i)_j}{(a+u_i)_j(b+u_{i+k})_j}
=\dfrac{(a+u_{i})_j}{(b+u_{i})_j}
\dfrac{(a+u_{i+j})_k(b+u_i)_k}{(b+u_{i+j})_k(a+u_i)_k}. 
\end{align}
The $(n+1)\times(n+1)$ minor determinants of the matrix 
\begin{align}
G=\Big(\dfrac{(c+u_i)_j}{(d+u_i)_j}\Big)_{0\le i\le N,\, 0\le j\le n}
\end{align}
can be computed by means of a special case of Krattenthaler's 
determinant formula \cite{K1999} (see Appendix B).  
In fact by \eqref{eq:B3}, we have 
\begin{align}
\det G^{i,i+1,\ldots,i+n}_{0,1,\ldots,n}&=
\dfrac{\tprod{l=1}{n}\,l!(d-c)_l}{\tprod{l=0}{n}(d+u_{i+l})_n},
\nonumber\\
\det G^{i,\ldots,\widehat{i+k},\ldots,i+n+1}_{0,1,\ldots,n}
&=\dfrac{\tprod{l=1}{n}\,l!(d-c)_l}{\tprod{l=1}{n+1}(d+u_{i+l})_n}\,
\,\dfrac{(-1)^k(-n-1)_k(d+u_{i+n})_k}{k!\,(d+u_i)_k}. 
\end{align}
Hence $U_{ij}$ is computed as 
\begin{align}
U_{ij}&=
\dfrac{(a+u_i)_j}{(b+u_i)_j}
\dsum{k=0}{n+1}
\dfrac{(-n-1)_k(d+u_{i+n})_k}{k!\ (d+u_i)_k}
\dfrac{(a+u_{i+j})_k(b+u_i)_k}{(a+u_{i})_k(b+u_{i+j})_k}
\dfrac{\mu_{i+k}\,\lambda_{i}}{\lambda_{i+k}\,\mu_i}. 
\end{align}  
The corresponding formula for $V_{ij}$ is obtained by exchanging the 
roles of $(m, n)$, $(a,c)$ and $(b,d)$.

\par\medskip
If we choose the prescribed values appropriately, the series $U_{ij}$ and 
$V_{ij}$ in Theorem \ref{thm:HG} give rise to {\em generalized hypergeometric 
series}
\begin{align}
{}_{r+1}F_r\left[\begin{matrix}
\alpha_0,\alpha_1,\ldots,\alpha_r\\
\beta_1,\ \ldots.\ \beta_r
\end{matrix};\ z\right]
=\dsum{k=0}{\infty}
\dfrac{(\alpha_0)_k(\alpha_1)_k\cdots(\alpha_r)_k}{(1)_k(\beta_1)_k\cdots(\beta_r)_k}
z^k.  
\end{align}
Consider the case where
\begin{align}
v_k=\dfrac{\lambda_k}{\mu_k}=
\dfrac{(s_1)_k\cdots(s_r)_k}{(t_1)_k\cdots(t_r)_k} \left(\dfrac{z}{w}\right)^k 
\qquad(k=0,1,2,\ldots,N)
\end{align}
with complex parameters $s_1,\ldots,s_r$ and $t_1,\ldots,t_r$. 
Since
\begin{align}
\dfrac
{\mu_{i+k}\lambda_i}
{\lambda_{i+k}\mu_i}=
\dfrac{(t_1+i)_k\cdots(t_r+i)_k}{(s_1+i)_k\cdots(s_r+i)_k} \left(\dfrac{w}{z}\right)^k, 
\end{align}
$U_{ij}$ and $V_{ij}$ are determined as 
\begin{align}
U_{ij}&=
\dfrac{(a+u_i)_j}{(b+u_i)_j}\,
{}_{r+4}F_{r+3}\left[\begin{matrix}
-n-1, d+u_{i+n}, a+u_{i+j}, b+u_i, t_1+i,\ldots,t_r+i\\
d+u_i, a+u_{i}, b+u_{i+j}, s_1+i,\ldots,s_r+i\\
\end{matrix}; \dfrac{w}{z}\right], 
\nonumber\\
V_{ij}&=
\dfrac{(c+u_i)_j}{(d+u_i)_j}\,
{}_{r+4}F_{r+3}\left[\begin{matrix}
-m-1, b+u_{i+n}, c+u_{i+j}, d+u_i, s_1+i,\ldots,s_r+i\\
b+u_i, c+u_{i}, d+u_{i+j}, t_1+i,\ldots,t_r+i\\
\end{matrix}; \dfrac{z}{w}\right]. 
\end{align}
If we choose the prescribed values
\begin{align}
v_k=\dfrac{\lambda_k}{\mu_k}=
\dfrac{(b+u)_k(c+u)_k}{(a+u)_k(d+u)_k}
\dfrac{(s_1)_k\cdots(s_r)_k}{(t_1)_k\cdots(t_r)_k} \left(\dfrac{z}{w}\right)^k 
\qquad(k=0,1,2,\ldots,N), 
\end{align}
then $U_{ij}$ and $V_{ij}$ are slightly simplified as
\begin{align}
U_{ij}&=
\dfrac{(a+u_i)_j}{(b+u_i)_j}\,
{}_{r+3}F_{r+2}\left[\begin{matrix}
-n-1, d+u_{i+n}, a+u_{i+j}, t_1+i,\ldots,t_r+i\\
c+u_i,, b+u_{i+j}, s_1+i,\ldots,s_r+i\\
\end{matrix}; \dfrac{w}{z}\right], 
\nonumber\\
V_{ij}&=
\dfrac{(c+u_i)_j}{(d+u_i)_j}\,
{}_{r+3}F_{r+2}\left[\begin{matrix}
-m-1, b+u_{i+n}, c+u_{i+j}, s_1+i,\ldots,s_r+i\\
a+u_i, d+u_{i+j}, t_1+i,\ldots,t_r+i\\
\end{matrix}; \dfrac{z}{w}\right]. 
\end{align}

\par\medskip
As for the Pad\'e interpolation problem for the rational functions 
$f_j(x)$ and $g_j(x)$ as in 
\eqref{eq:funcfg}, one can construct another type of determinant formula 
for $P_m(x)$ and $Q_n(x)$ involving hypergeometric series.  

\begin{thm}\label{thm:HG2}
Consider the Pad\'e interpolation problem \eqref{eq:PadeR0}, \eqref{eq:PadeR1}
for the rational functions 
$f_j(x)$, $g_j(z)$ in \eqref{eq:funcfg} 
and the reference points $u_k=u+k$ $(k=0,1,\ldots,N; N=m+n)$. 
Then the solution $(P_m(x),Q_n(x))$ of Theorem 
\ref{thm:PQGen} is expressed as 
\begin{align}\label{eq:PDetHG2}
\arraycolsep=2pt
P_m(x)&=
K_{m,n}(c,d)
\dprod{i=0}{N}\,\lambda_{i}\,
\det\left[\,
\begin{matrix}
f_0(x) & \ldots & f_m(x) \\[4pt]
\Phi_{0,0} & \ldots & \Phi_{0,m} \\
\vdots &&\vdots\\
\Phi_{m-1,0}&\ldots & \Phi_{m-1,m}
\end{matrix}\,
\right], 
\\
\label{eq:QDetHG2}
Q_n(x)&=
\epsilon_{m,n}
K_{n,m}(a,b)
\dprod{i=0}{N}\,\mu_{i}\,
\det\left[\,
\begin{matrix}
g_0(x) & \ldots & g_n(x) \\[4pt]
\Psi_{0,0} & \ldots & \Psi_{0,n} \\
\vdots &&\vdots\\
\Psi_{n-1,0}&\ldots & \Psi_{n-1,n}
\end{matrix}\,
\right], 
\end{align}
where 
\begin{align}
\Phi_{ij}&=
\dfrac{(a+u)_j}{(b+u)_j}\dsum{k=0}{N}
\dfrac{(-N)_k(d+u+N-1-i)_k}
{(1)_k(c+u-i)_k}
\dfrac{(c+u)_k}{(d+u)_k}
\dfrac{(a+u_j)_k(b+u)_k}{(b+u_j)_k(a+u)_k}
\dfrac{\mu_k}{\lambda_k}, 
\\
\Psi_{ij}&=
\dfrac{(c+u)_j}{(d+u)_j}\dsum{k=0}{N}
\dfrac{(-N)_k(b+u+N-1-i)_k}
{(1)_k(a+u-i)_k}
\dfrac{(a+u)_k}{(b+u)_k}
\dfrac{(c+u_j)_k(d+u)_k}{(d+u_j)_k(c+u)_k}
\dfrac{\lambda_k}{\mu_k}.  
\end{align}  
\end{thm}

We remark that if the prescribed values are given by 
\begin{align}
v_k=\dfrac{\lambda_k}{\mu_k}=
\dfrac{(s_1)_k\cdots(s_r)_k}{(t_1)_k\cdots(t_r)_k} \left(\dfrac{z}{w}\right)^k 
\qquad(k=0,1,2,\ldots,N), 
\end{align}
then 
$\Phi_{ij}$ and $\Psi_{ij}$ give rise to generalized hypergeometric series  
\begin{align}
\Phi_{ij}&=
\dfrac{(a+u)_j}{(b+u)_j}\,
{}_{r+5}F_{r+4}\left[\begin{matrix}
-N, d+u+N-1-i, c+u, a+u+j, b+u, t_1,\ldots,t_r\\
c+u-i, d+u, b+u+j, a+u, s_1,\ldots,s_r\\
\end{matrix}; \dfrac{w}{z}\right], 
\nonumber\\
\Psi_{ij}&=
\dfrac{(c+u)_j}{(d+u)_j}\,
{}_{r+5}F_{r+4}\left[\begin{matrix}
-N, b+u+N-1-i, a+u, c+u+j, d+u, s_1,\ldots,s_r\\
a+u-i, b+u, d+u+j, c+u, t_1,\ldots,t_r\\
\end{matrix}; \dfrac{z}{w}\right].  
\end{align}
If we choose 
\begin{align}
v_k=\dfrac{\lambda_k}{\mu_k}=
\dfrac{(b+u)_k(c+u)_k}{(a+u)_k(d+u)_k}
\dfrac{(s_1)_k\cdots(s_r)_k}{(t_1)_k\cdots(t_r)_k} \left(\dfrac{z}{w}\right)^k 
\qquad(k=0,1,2,\ldots,N), 
\end{align}
then $\Phi_{ij}$ and $\Psi_{ij}$ are simplified as
\begin{align}
\Phi_{ij}&=
\dfrac{(a+u)_j}{(b+u)_j}\,
{}_{r+3}F_{r+2}\left[\begin{matrix}
-N, d+u+N-1-i, a+u+j, t_1,\ldots,t_r\\
c+u-i, b+u+j, s_1,\ldots,s_r\\
\end{matrix}; \dfrac{w}{z}\right], 
\nonumber\\
\Psi_{ij}&=
\dfrac{(c+u)_j}{(d+u)_j}\,
{}_{r+3}F_{r+2}\left[\begin{matrix}
-N, b+u+N-1-i, c+u+j, s_1,\ldots,s_r\\
a+u-i, d+u+j,  t_1,\ldots,t_r\\
\end{matrix}; \dfrac{z}{w}\right].  
\end{align}

\par\medskip
In order to obtain the expression of Theorem \ref{thm:HG2}, 
we first rewrite \eqref{eq:PDet} as 
\begin{align}
\arraycolsep=2pt
P_m(x)
&=\dprod{i=0}{N}\lambda_i\ 
\det
\left[\,
\begin{matrix}
f_0(x) & \cdots & f_m(x) &
0 & \cdots & 0 \\[4pt]
\dfrac{\mu_0}{\lambda_0} f_0(u_0) & \cdots & \dfrac{\mu_0}{\lambda_0} f_m(u_0) &
g_0(u_0) & \cdots & g_n(u_0) \\
\vdots &&\vdots&\vdots&&\vdots\\
\dfrac{\mu_N}{\lambda_N} f_0(u_N) & \cdots & 
\dfrac{\mu_N}{\lambda_N} f_m(u_N) &
g_0(u_N) & \cdots & g_n(u_N)
\end{matrix}\,
\right]
\nonumber\\[8pt]
&=\dprod{i=0}{N}\lambda_i\ 
\det\left[
\begin{matrix}
f(x) & 0 \\[4pt]
\widetilde{F}  & G
\end{matrix}
\right]. 
\end{align}
We construct an $(N+1)\times(N+1)$ invertible 
matrix $L=\big(L_{ij}\big)_{i,j=0}^N$ 
such that $(LG)_{ij}=0$ for $i+j<N$, and 
define 
$M=\big(M_{ij}\big)_{i,j=0}^{n}$ by $M_{ij}=(LG)_{m+i,j}$.  
If we set $\Phi=L \widetilde{F}$, we have 
\begin{align}
\arraycolsep=2pt
\left[\begin{matrix}
\,1 & \\[6pt]
 & L\ 
\end{matrix}\right]
\left[\begin{matrix}
\,f(x) & 0\\[6pt]
\widetilde{F} & G \ 
\end{matrix}\right]=
\left[\begin{matrix}
\,f(x) & 0\\[6pt]
\Phi & LG\ 
\end{matrix}\right]
=
\left[\begin{matrix}
\,f(x) & \ 0\ \\[2pt]
\Phi' & 0\\
\Phi'' & M
\end{matrix}\right]. 
\end{align}
Hence, by taking the determinants of the both sides we obtain 
\begin{align}
P_m(x) =\lambda_0\cdots\lambda_N 
\det\left[\begin{matrix}
f_0(x)& \ldots & f_m(x)\\
\Phi_{00}& \ldots & \Phi_{0m}\\
\vdots & \ddots & \vdots\\
\Phi_{m-1,0}& \ldots & \Phi_{m-1,m}
\end{matrix}
\right]\dfrac{\det M}{\det L}, 
\end{align}
which will give formula \eqref{eq:PDetHG2} with 
$K_{m,n}(c,d)=\det M/\det L$. 
In view of 
\begin{align}
g_j(u_k)=
\dfrac{(c+u_k)_j}{(d+u_k)_j}
=\dfrac{(d+u)_k}{(c+u)_k}\dfrac{(c+u_j)_k}{(d+u_j)_k}
\dfrac{(c+u)_j}{(d+u)_j}, 
\end{align}
we recall the Saalsch\"utz sum
\begin{align}
{}_{3}F_{2}\left[\begin{matrix}
-N, d+u+N-1-i, c+u+j \\ 
c+u-i, d+u+j
\end{matrix};\,1\right]=
\dfrac{(d-c)_N(-i-j)_N}{(c+u-i)_N(d+u+j)_N}, 
\end{align}
namely, 
\begin{align}
\dsum{k=0}{N}\dfrac{(-N)_k(d+u+N-1-i)_k(c+u_j)_k}
{(1)_k(c+u-i)_k(d+u_j)_k}=
\dfrac{(d-c)_N(-i-j)_N}{(c+u-i)_N(d+u+j)_N}. 
\end{align}
With this observation, 
we define the matrix $L=\big(L_{ij}\big)_{i,j=0}^{N}$ 
by  
\begin{align}
L_{ij}=
\dfrac{(-N)_j(d+u+N-1-i)_j}
{(1)_j(c+u-i)_j}
\dfrac{(c+u)_j}{(d+u)_j}\qquad(0\le i,j\le N). 
\end{align}
The we have 
\begin{align}
(LG)_{ij}&=\dsum{k=0}{N}
\dfrac{(-N)_k(d+u+N-1-i)_k}
{(1)_j(c+u-i)_k}
\dfrac{(c+u_j)_k}{(d+u_j)_k}
\dfrac{(c+u)_j}{(d+u)_j}
\nonumber\\
&=
\dfrac{(d-c)_N(-i-j)_N}{(c+u-i)_N(d+u+j)_N}
\dfrac{(c+u)_j}{(d+u)_j}
\end{align}
by the Saalsch\"utz sum. 
In particular $(LG)_{ij}=0\quad(i+j<N)$.  The determinant of the matrix $M$ is 
computed as
\begin{align}
\det M=(-1)^{\binom{n+1}{2}} \dprod{j=0}{n}(LG)_{N-j,j}
=(-1)^{\binom{n+1}{2}}\dprod{j=0}{n}
\dfrac{(d-c)_N(-N)_N}{(c+u-N+j)_N(d+u+j)_N}
\dfrac{(c+u)_j}{(d+u)_j}. 
\end{align}
Also, the entires of $\Phi=L \widetilde{F}$ are expressed as 
\begin{align}
\Phi_{ij}&=\dsum{k=0}{N}L_{ik}\frac{\mu_k}{\lambda_k}f_j(u_k)
\nonumber\\
&=\dsum{k=0}{N}
\dfrac{(-N)_k(d+u+N-1-i)_k}
{(1)_k(c+u-i)_k}
\dfrac{(c+u)_k}{(d+u)_k}
\dfrac{\mu_k}{\lambda_k}\dfrac{(a+u_k)_j}{(b+u_k)_j}
\nonumber\\
&=
\dfrac{(a+u)_j}{(b+u)_j}\dsum{k=0}{N}
\dfrac{(-N)_k(d+u+N-1-i)_k}
{(1)_k(c+u-i)_k}
\dfrac{(c+u)_k}{(d+u)_k}
\dfrac{(a+u_j)_k(b+u)_k}{(b+u_j)_k(a+u)_k}
\dfrac{\mu_k}{\lambda_k}. 
\end{align}
The determinant of $L$ can be computed again by Krattenthaler's formula:
\begin{align}
\det L
&=
\det\left(
\dfrac{(d+u+N-1-i)_j}{(c+u-i)_j}
\right)_{i,j=0}^{N}
\dprod{j=0}{N}
\dfrac{(-N)_j}{(1)_j}\dfrac{(c+u)_j}{(d+u)_j}
\nonumber\\
&=
(-1)^{\binom{N+1}{2}}
\dprod{j=0}{N}
\dfrac{(c-d-N+1)_j}{(c+u-j)_N}
\dfrac{(-N)_j(c+u)_j}{(d+u)_j}. 
\end{align}
The constant factor in \eqref{eq:PDetHG2} is determined as 
$K_{m,n}(c,d)=\det M/\det L$.  
%%%%%%%%%%%%%%%%%%%%%%%%%%%%%%%%%%%%%%%%%%%%%%

%%%%%%%%%%%%%%%%%%%%%%%%%%%%%%%%%%
\section{Three types of very well-poised hypergeometric series}
%%%%%%%%%%%%%%%%%%%%%%%%%%%%%%%%%%
In this section we consider three classes of hypergeometric series
\par\medskip
\centerline{
\begin{tabular}{lcl}
(0) \ rational & $\ldots$ & ordinary hypergeometric series\\
(1) \ trigonometric & $\ldots$ & basic (or $q$-)hypergeometric series\\
(2) \ elliptic & $\ldots$ & elliptic hypergeometric series
\end{tabular}
}
\par\medskip\noindent
corresponding to the choice of a ``fundamental'' function $[x]$: 
\par\medskip
\centerline{
\begin{tabular}{lcll}
(0) \ rational & $:$ & $[x]=e^{c_0x^2+c_1}\, x$ & ($\Omega=0$)\\
(1) \ trigonometric & $:$ & $[x]=e^{c_0x^2+c_1}\,\sin(\pi x/\omega)$
& ($\Omega=\mathbb{Z}\omega$)\\
(2) \ elliptic & $:$ & $[x]=e^{c_0x^2+c_1}\,\,\sigma(x|\Omega)$ & 
($\Omega=
\mathbb{Z}\omega_1\oplus \mathbb{Z}\omega_2$)  
\end{tabular}
}
\par\medskip\noindent
where $\sigma(x|\Omega)$ is the Weierstrass sigma function associated 
with the period lattice 
$\Omega=\mathbb{Z}\omega_1\oplus \mathbb{Z}\omega_2$. 
It is known that 
these classes of functions $[x]$ are characterized 
by the so-called {\em Riemann relation}: 
For any $x,\alpha,\beta,\gamma\in\mathbb{C}$, 
\begin{equation}
[x+\alpha][x-\alpha][\beta+\gamma][\beta-\gamma]
+
[x+\beta][x-\beta][\gamma+\alpha][\gamma-\alpha]
+
[x+\gamma][x-\gamma][\alpha+\beta][\alpha-\beta]
=0. 
\end{equation}
By the notation $[x\pm y]=[x+y][x-y]$ of the product 
of two factors, this relation is expressed as
\begin{equation}
[x\pm\alpha][\beta\pm\gamma]
+
[x\pm\beta][\gamma\pm\alpha]
+
[x\pm\gamma][\alpha\pm\beta]
=0. 
\end{equation}
In what follows,  we fix 
a nonzero entire function $[x]$ satisfying this functional equation. 

Fixing a generic constant $\delta$, 
we define the $\delta$-{\em shifted factorials}  
$[x]_k$ by 
\begin{equation}
[x]_k=[x]_{\delta,k}=[x][x+\delta]\cdots [x+(k-1)\delta]\qquad(k=0,1,2,\ldots).  
\end{equation}
Then we define the {\em very well-poised} hypergeometric series 
${}_{r+5}V_{r+4}\big[
a_0; a_1\, \cdots\, a_r 
\big| z
\big]$ 
associated with $[x]$ by 
\begin{equation}
{}_{r+5}V_{r+4}\Big[
a_0; a_1,\,\cdots,a_r\, 
\Big|\,z
\Big]
=
\dsum{k=0}{\infty}
\dfrac{[a_0+2k\delta]}{[a_0]}
\dfrac{[a_0]_k\ [a_1]_k\ \cdots\ [a_r]_k}{[\delta]_k[\delta+a_0-a_1]_k
\cdots[\delta+a_0-a_r]_k}\,z^k.  
\end{equation}
In this paper we use this notation only for terminating series assuming 
that $a_i$ is of the form $-n\delta$ ($n=0,1,2,\ldots$) for some $i$. 
When $z=1$ we also write 
\begin{equation}
{}_{r+5}V_{r+4}\Big[
a_0; a_1,\, \cdots, a_r\Big]
=
\dsum{k=0}{\infty}
\dfrac{[a_0+2k\delta]}{[a_0]}
\dfrac{[a_0]_k\ [a_1]_k\ \cdots\ [a_r]_k}{[\delta]_k[\delta+a_0-a_1]_k
\cdots[\delta+a_0-a_r]_k}. 
\end{equation}
In this notation, the celebrated Frenkel-Turaev sum is 
expressed as 
\begin{align}
&{}_{10}V_{9}\Big[a_0;a_1,a_2,a_3,a_4,a_5\Big]
\nonumber\\
&
=
\dfrac{
[\delta+a_0]_N[\delta+a_0-a_1-a_2]_N[\delta+a_0-a_1-a_3]_N[\delta+a_0-a_2-a_3]_N
}{
[\delta+a_0-a_1]_N
[\delta+a_0-a_2]_N
[\delta+a_0-a_3]_N
[\delta+a_0-a_1-a_2-a_3]_N
}, 
\end{align}
under the balancing condition $a_1+\cdots+a_5=2a_0+\delta$ and 
the termination condition $a_5=-N\delta$ ($N=0,1,2,\ldots$). 
(See for example \cite{GR2004}, \cite{KN2003}.)  

We remark that, in the rational case where $[x]=x$ and $\delta=1$, 
the ${}_{r+5}V_{r+4}$ series defined above is expressed in terms of a 
${}_{r+2}F_{r+1}$-series: 
\begin{align}
{}_{r+5}V_{r+4}\Big[
a_0; a_1,\,\cdots,a_r\, 
\Big|\,z
\Big]
={}_{r+2}F_{r+1}\left[\begin{matrix}
a_0, \tfrac{1}{2}a_0+1, a_1, \ldots a_r\\[2pt]
\tfrac{1}{2}a_0, \ b_1,\ \ldots, \ b_r
\end{matrix};\,z \right]
\end{align}
where $b_i=1+a_0-a_i$ $(i=1,\ldots,r)$.  
Also, in the trigonometric case where $[x]=e^{cx/2}-e^{-cx/2}$, 
\begin{align}
{}_{r+5}V_{r+4}\Big[
a_0; a_1,\,\cdots,a_r\, 
\Big|\,z
\Big]
&=
\dsum{k=0}{\infty}
\dfrac{1-q^{2k}t_0}{1-t_0}
\dfrac{(t_0;q)_k\ (t_1;q)_k\ \cdots\ (t_r;q)_k}
{(q;q)_k(qt_0/t_1;q)_k\cdots(qt_0/t_r)_k}\,s^k
\nonumber\\
&=
{}_{r+3}W_{r+2}\Big[t_0; t_1,\ldots,t_r;q,s\Big]
\end{align}
in the notation of very well-poised $q$-hypergeometric series \cite{GR2004}, 
where $q=e^{c\delta}$, $t_i=e^{cx_i}$ ($i=0,1,\ldots,r$) and 
$s=(qt_0)^{\frac{r-1}{2}}z/t_1\cdots t_r$.  
We discuss below a class of Pad\'e interpolation problems that can be 
formulated in an unified manner in the three types of very well-poised 
hypergeometric series.  

\par\medskip
Taking the two sequence of meromorphic functions 
\begin{align}\label{eq:Efg}
f_j(x)&=
\dfrac{[a\pm x]_j}{[b\pm x]_j}
=
\dfrac{[a+x]_j[a-x]_j}{[b+x]_j[b-x]_j}
,\quad
\nonumber\\
g_j(x)&=\dfrac{[c\pm x]_j}{[d\pm x]_j}
=
\dfrac{[c+x]_j[c-x]_j}{[d+x]_j[d-x]_j}
\quad(j=0,1,2,\ldots)
\end{align}
and the reference points $u_k=u+k\delta$ ($k=0,1,2,\ldots$), 
we consider the Pad\'e interpolation problem 
\begin{align}
\label{eq:PadeE0}
&&\dfrac{P_m(u_k)}{Q_n(u_k)}=v_k=\dfrac{\lambda_k}{\mu_k}
\qquad(k=0,1,\ldots,N)
\end{align}
for a pair of functions 
\begin{align}
\label{eq:PadeE1}
P_m(x)&=p_{m,0}\,f_0(x)+p_{m,1}\,f_1(x)+\cdots+p_{m,m}\,f_m(x),
\nonumber\\
Q_n(x)&=q_{n,0}\,g_0(x)+q_{n,1}\,g_1(x)+\cdots+q_{n,n}\,g_n(x)
%P_m(u)=\dsum{j=0}{m}\,p_{m,j}\,f_j(u),
%\quad
%Q_n(u)=\dsum{j=0}{n}\,q_{n,j}\,g_j(u). 
\end{align}
where $N=m+n$.  
The prescribed values $v_k=\lambda_k/\mu_k$
($k=0,1,2\ldots,N$) will be specified later. 
%v_k=z^{k}
%\dfrac{[\delta\!-\!a\!+\!u_0]_k[b+u_0]_k}{[a+u_0]_k[\delta\!-\!b\!+\!u_0]_k}
%\dfrac{[c+u_0]_k[\delta\!-\!d+u_0]_k}{[\delta\!-\!c\!+\!u_0]_k[d+u_0]_k}\,
%\dprod{s=1}{r}\dfrac{[\delta\!-\!e_s\!+\!u_0]_k}{[e_s+u_0]_k}. 

\begin{thm}\label{thm:EHG}
Consider the Pad\'e interpolation problem \eqref{eq:PadeE0}, \eqref{eq:PadeE1}
for the functions 
$f_j(x)$, $g_j(x)$ in \eqref{eq:Efg} 
and the reference points $u_k=u+k\delta$ $(k=0,1,\ldots,N; N=m+n)$. 
Then the solution $(P_m(x),Q_n(x))$ of Theorem 
\ref{thm:PQGen} is explicitly given by 
\begin{align}\label{eq:PDetUEHG}
\arraycolsep=2pt
P_m(x)&=
C_n(c,d)
\dfrac{\tprod{l=1}{n}[2u_m+l\delta]_l[\delta]_l}
{\tprod{l=0}{n}\,[d\pm u_{m+l}]_n}\,
\dprod{i=0}{m-1}\,\mu_i\,\dprod{i=0}{n}\,\lambda_{m+i}\,
\det\left[\,
\begin{matrix}
f_0(x) & \ldots & f_m(x) \\[4pt]
U_{0,0} & \ldots & U_{0,m} \\
\vdots &&\vdots\\
U_{m-1,0}&\ldots & U_{m-1,m}
\end{matrix}\,
\right], 
\\
\label{eq:QDetVEHG}
Q_n(x)&=
\epsilon_{m,n}
C_m(a,b)
\dfrac{\tprod{l=1}{m}[2u_n+l\delta]_l[\delta]_l}
{\tprod{l=0}{m}\,[b\pm u_{n+l}]_m}\,
\dprod{i=0}{n-1}\lambda_i\dprod{i=0}{m}\mu_{n+i}
\det\left[\,
\begin{matrix}
g_0(x) & \ldots & g_n(x) \\[4pt]
V_{0,0} & \ldots & V_{0,n} \\
\vdots &&\vdots\\
V_{n-1,0}&\ldots & U_{n-1,n}
\end{matrix}\,
\right], 
\end{align}
where 
\begin{align}
C_n(c,d)=(-1)^{\binom{n+1}{2}}\dprod{k=1}{n} [d-c]_k[c+d+(k-1)\delta]_k,
\quad \epsilon_{m,n}=(-1)^{mn+m+n}
\end{align}
and 
\begin{align}
U_{ij}&=
\dfrac{[a\pm u_{i}]_{j}}{[b\pm u_{i}]_{j}}
\dsum{k=0}{n+1}
\dfrac{[2u_i+2k\delta]}{[2u_i]}
\dfrac{[2u_i]_k\ [-(n+1)\delta]_k}{[\delta]_k[2u_i+(n+2)\delta]_k}
\dfrac{[u_i-d+\delta]_k\ [u_i+d+n\delta]_k}{[u_i+d]_k[u_i-d+(1-n)\delta]_k}
\nonumber\\
&\qquad\qquad\qquad\cdot
\dfrac{[u_i-a+\delta]_k[u_{i}+a+j\delta]_{k}[u_i+b]_k[u_i-b+(1-j)\delta]_k}
{[u_i+a]_k[[u_i-a+(1-j)\delta]_k[u_i-b+\delta]_k[u_{i}+b+j\delta]_k}
\dfrac{\mu_{i+k}\lambda_i}{\lambda_{i+k}\mu_i}, 
\\[8pt]
V_{ij}&=
\dfrac{[c\pm u_{i}]_{j}}{[d\pm u_{i}]_{j}}
\dsum{k=0}{m+1}
\dfrac{[2u_i+2k\delta]}{[2u_i]}
\dfrac{[2u_i]_k\ [-(m+1)\delta]_k}{[\delta]_k[2u_i+(m+2)\delta]_k}
\dfrac{[u_i-b+\delta]_k\ [u_i+b+m\delta]_k}{[u_i+b]_k[u_i-b+(1-m)\delta]_k}
\nonumber\\
&\qquad\qquad\qquad\cdot
\dfrac{[u_i-c+\delta]_k[u_{i}+c+j\delta]_{k}[u_i+d]_k[u_i-d+(1-j)\delta]_k}
{[u_i+c]_k[[u_i-c+(1-j)\delta]_k[u_i-d+\delta]_k[u_{i}+d+j\delta]_k}
\dfrac{\lambda_{i+k}\mu_i}{\mu_{i+k}\lambda_i}. 
\end{align}
\end{thm}

\par\medskip
As before we consider the expansion 
\begin{align}
U_{ij}&
=\dsum{k=0}{n+1}(-1)^k 
\,\dfrac{G^{i,\ldots,\widehat{i+k},\ldots,i+n+1}_{0,1,\ldots,n}}
{G^{i+1,\ldots,i+n+1}_{0,1,\ldots,n}}\,
f_j(u_{i+k})
\dfrac{\mu_{i+k}\,\lambda_{i}}{\lambda_{i+k}\,\mu_{i}} 
\end{align}
of the determinant of \eqref{eq:defU}. 
In this case we have
\begin{align}
f_j(u_{i+k})=
\dfrac{[a\pm u_{i+k}]_j}{[b\pm u_{i+k}]_{j}}
=
\dfrac{[a\pm u_{i}]_{j}}{[b\pm u_{i}]_{j}}
\dfrac{[u_i-a+\delta]_k[u_{i}+a+j\delta]_{k}[u_i+b]_k[u_i-b+(1-j)\delta]_k}
{[u_i+a]_k[[u_i-a+(1-j)\delta]_k[u_i-b+\delta]_k[u_{i}+b+j\delta]_k}. 
\end{align}
The $(n+1)\times (n+1)$ minor determinants of 
the matrix
\begin{align}
G=
\big(g_j(u_i)\big)_{0\le i\le N,\, 0\le j\le n}
=
\left(
\dfrac{[c\pm u_i]_j}{[d\pm u_i]_j}
\right)_{0\le i\le N,\, 0\le j\le n} 
\end{align}
can be computed by means of an elliptic extension of 
Krattenthaler's formula (see Appendix B). 
In fact, by \eqref{eq:B14} we have
\begin{align}
\det G^{i,i+1,\ldots,i+n}_{0,1,\ldots,n}
&=
C_n(c,d)
\dfrac{\tprod{l=1}{n}[2u_i+l\delta]_l[\delta]_l}
{\tprod{l=0}{n}\,[d\pm u_{i+l}]_n}, 
\nonumber\\
%&C_n(c,d)=(-1)^{\binom{n+1}{2}}\dprod{k=1}{n} [d-c]_k[c+d+(k-1)\delta]_k, 
%\end{align}
%and
%\begin{align}
\det G^{i,\ldots,\widehat{i+k},\ldots ,i+n+1}_{0,1,\ldots,n}
&=\,
C_n(c,d)\,
\dfrac{\tprod{l=1}{n}\,[2u_{i+1}+l\delta]_l[\delta]_l}
{\tprod{l=0}{n}[d\pm u_{i+1+l}]_n}
\nonumber\\[4pt]
&\quad \cdot (-1)^k\,
\dfrac{[2u_i+2k\delta]}{[2u_i]}
\dfrac{[2u_i]_k[-(n+1)\delta]_k}{[\delta]_k[2u_i+(n+2)\delta]_k}
\dfrac{[u_i-d+\delta]_k[u_i+d+n\delta]_k}{[u_i+d]_k[u_i-d+(1-n)\delta]_k}. 
\end{align}
Hence $U_{ij}$ is computed as 
\begin{align}
U_{ij}&=
\dfrac{[a\pm u_{i}]_{j}}{[b\pm u_{i}]_{j}}
\dsum{k=0}{n+1}
\dfrac{[2u_i+2k\delta]}{[2u_i]}
\dfrac{[2u_i]_k\ [-(n+1)\delta]_k}{[\delta]_k[2u_i+(n+2)\delta]_k}
\dfrac{[u_i-d+\delta]_k\ [u_i+d+n\delta]_k}{[u_i+d]_k[u_i-d+(1-n)\delta]_k}
\nonumber\\
&\qquad\qquad\quad\cdot
\dfrac{[u_i-a+\delta]_k[u_{i}+a+j\delta]_{k}[u_i+b]_k[u_i-b+(1-j)\delta]_k}
{[u_i+a]_k[[u_i-a+(1-j)\delta]_k[u_i-b+\delta]_k[u_{i}+b+j\delta]_k}
\dfrac{\mu_{i+k}\lambda_i}{\lambda_{i+k}\mu_i}. 
\end{align}
The corresponding formula for $V_{ij}$ is obtained by exchanging 
the roles of $(m,n)$, $(a,c)$ and $(b,d)$. 
\par\medskip
Consider the case where the prescribed values are specified as 
\begin{align}
v_k=\dfrac{\lambda_k}{\mu_k}=
\Big(\dfrac{z}{w}\Big)^k
\dprod{s=1}{r}
\dfrac{[u-e_s+\delta]_k}{[u+e_s]_k}\qquad(k=0,1,\ldots,N).  
\end{align}
Then we obtain very well-poised series
\begin{align}
U_{ij}&=\dfrac{[a\pm u_i]_j}{[b\pm u_i]_j}
\,{}_{r+12}V_{r+11}\Big[2u_i; -(n+1)\delta, u_i-d+\delta, u_i+d+n\delta, 
\nonumber\\
&\qquad\qquad\qquad
u_i-a+\delta, u_i+a+j\delta, 
u_i+b, u_i-b+(1-j)\delta, u_i+e_1,\ldots,u_i+e_r
\Big|\dfrac{w}{z}\Big],
\\
V_{ij}&=\dfrac{[c\pm u_i]_j}{[d\pm u_i]_j}
\,{}_{r+12}V_{r+11}\Big[2u_i; -(m+1)\delta, u_i-b+\delta, u_i+b+n\delta, 
\nonumber\\
&\qquad\qquad\qquad
u_i\!-c\!+\!\delta, u_i\!+\!c\!+\!j\delta, 
u_i\!+d\!, u_i\!-\!d\!+\!(1-j)\delta, u_i\!-\!e_1\!+\!\delta,\ldots,u_i\!-\!e_r\!+\!\delta
\Big|\dfrac{z}{w}\Big]. 
\end{align}
When 
\begin{align}
v_k=\dfrac{\lambda_k}{\mu_k}=
\Big(\dfrac{z}{w}\Big)^k
\dfrac{[u-a+\delta]_k[u+b]_k}{[u+a]_k[u-b+\delta]_k}
\dfrac{[u+c]_k[u-d+\delta]_k}{[u-c+\delta]_k[u+d]_k}
\dprod{s=1}{r}
\dfrac{[u-e_s+\delta]_k}{[u+e_s]_k}
\end{align}
we obtain simpler very well-poised hypergeometric series
\begin{align}
U_{ij}&=\dfrac{[a\pm u_i]_j}{[b\pm u_i]_j}
\,{}_{r+10}V_{r+9}\Big[2u_i; -(n+1)\delta, u_i-c+\delta, u_i+d+n\delta, 
\nonumber\\
&\qquad\qquad\qquad\quad
u_i+a+j\delta, u_i-b+(1-j)\delta, 
u_i+e_1,\ldots,u_i+e_r
\Big|\dfrac{w}{z}\Big],
\nonumber\\
V_{ij}&=\dfrac{[c\pm u_i]_j}{[d\pm u_i]_j}
\,{}_{r+10}V_{r+9}\Big[2u_i; -(m+1)\delta, u_i-a+\delta, u_i+b+n\delta, 
\nonumber\\
&\qquad\qquad\qquad\quad
u_i+c+j\delta, u_i-d+(1-j)\delta, 
u_i-e_1+\delta,\ldots,u_i-e_r+\delta
\Big|\dfrac{z}{w}\Big]. 
\end{align}

\par\medskip
Another type of determinantal expression
for $P_m(x)$ and $Q_n(x)$ is formulated as follows. 
We remark that this type of determinant formulas has also been 
discussed in \cite{NTY2013}.  
In what follows, we use the notation
\begin{align}
V^{(k)}\Big[a_0; a_1,\, \cdots, a_r\Big]
=
\dfrac{[a_0+2k\delta]}{[a_0]}
\dfrac{[a_0]_k\ [a_1]_k\ \cdots\ [a_r]_k}{[\delta]_k[\delta+a_0-a_1]_k
\cdots[\delta+a_0-a_r]_k}. 
\end{align}
\begin{thm}\label{thm:EHG2}
Consider the Pad\'e interpolation problem \eqref{eq:PadeE0}, \eqref{eq:PadeE1}
for the functions 
$f_j(x)$, $g_j(z)$ in \eqref{eq:Efg} 
and the reference points $u_k=u+k\delta$ $(k=0,1,\ldots,N; N=m+n)$. 
Then the solution $(P_m(x),Q_n(x))$ of Theorem 
\ref{thm:PQGen} is expressed as 
\begin{align}\label{eq:PDetEHG2}
\arraycolsep=2pt
P_m(x)&=
K_{m,n}(c,d)
\dprod{i=0}{N}\,\lambda_{i}\,
\det\left[\,
\begin{matrix}
f_0(x) & \ldots & f_m(x) \\[4pt]
\Phi_{0,0} & \ldots & \Phi_{0,m} \\
\vdots &&\vdots\\
\Phi_{m-1,0}&\ldots & \Phi_{m-1,m}
\end{matrix}\,
\right], 
\\
\label{eq:QDetEHG2}
Q_n(x)&=
\epsilon_{m,n}
K_{n,m}(a,b)
\dprod{i=0}{N}\,\mu_{i}\,
\det\left[\,
\begin{matrix}
g_0(x) & \ldots & g_n(x) \\[4pt]
\Psi_{0,0} & \ldots & \Psi_{0,n} \\
\vdots &&\vdots\\
\Psi_{n-1,0}&\ldots & \Psi_{n-1,n}
\end{matrix}\,
\right], 
\end{align}
where 
\begin{align}
\Phi_{ij}&=
\dfrac{[a\pm u]_j}{[b\pm u]_j}\dsum{k=0}{N}\,
V^{(k)}\Big[2u; -N\delta, u-c+\delta+i\delta,u+d+(N-1)\delta-i\delta, u+c, u-d+\delta,
\nonumber\\
&\qquad \qquad\qquad\qquad\qquad
u+a+j\delta,u-b+\delta-j\delta,u-a+\delta,u+b\Big]\, 
\dfrac{\mu_k}{\lambda_k}, 
\\
\Psi_{ij}&=
\dfrac{[c\pm u]_j}{[d\pm u]_j}\dsum{k=0}{N}\,
V^{(k)}\Big[2u; -N\delta, u-a+\delta+i\delta,u+b+(N-1)\delta-i\delta, u+a, u-b+\delta,
\nonumber\\
&\qquad \qquad\qquad\qquad\qquad
u+c+j\delta,u-d+\delta-j\delta,u-c\delta,u+d\Big]\, 
\dfrac{\lambda_k}{\mu_k}. 
\end{align}  
\end{thm}

Consider the case where the prescribed values are specified as 
\begin{align}
v_k=\dfrac{\lambda_k}{\mu_k}=
\Big(\dfrac{z}{w}\Big)^k
\dprod{s=1}{r}
\dfrac{[u-e_s+\delta]_k}{[u+e_s]_k}\qquad(k=0,1,\ldots,N).  
\end{align}
Then we obtain very well-poised series
\begin{align}
\Phi_{ij}&=
\dfrac{[a\pm u]_j}{[b\pm u]_j}
{}_{r+14}V_{r+13}\Big[2u; -N\delta, u-c+\delta+i\delta,
u+d+(N-1)\delta-i\delta, u+c, u-d+\delta,
\nonumber\\
&\qquad \qquad\qquad\qquad
u+a+j\delta,u-b+\delta-j\delta,u-a+\delta,u+b,u+e_1,\ldots,u+e_r\Big|\dfrac{w}{z}\Big], 
\\
\Psi_{ij}&=
\dfrac{[c\pm u]_j}{[d\pm u]_j}
{}_{r+14}V_{r+13}\Big[2u; -N\delta, u-a+\delta+i\delta,
u+b+(N-1)\delta-i\delta, u+a, u-b+\delta,
\nonumber\\
&\qquad \qquad\qquad\qquad
u\!+\!c\!+\!j\delta,u\!-\!d\!+\!\delta\!-\!j\delta,
u\!-\!c\!+\!\delta,u\!+\!d,u\!-\!e_1\!+\!\delta,\ldots,u\!-\!e_r\!+\!\delta
\Big|\dfrac{z}{w}\Big]. 
\end{align}  
When 
\begin{align}
v_k=\dfrac{\lambda_k}{\mu_k}=
\Big(\dfrac{z}{w}\Big)^k
\dfrac{[u-a+\delta]_k[u+b]_k}{[u+a]_k[u-b+\delta]_k}
\dfrac{[u+c]_k[u-d+\delta]_k}{[u-c+\delta]_k[u+d]_k}
\dprod{s=1}{r}
\dfrac{[u-e_s+\delta]_k}{[u+e_s]_k}
\end{align}
we obtain simpler very well-poised hypergeometric series
\begin{align}
\Phi_{ij}&=
\dfrac{[a\pm u]_j}{[b\pm u]_j}
{}_{r+10}V_{r+9}\Big[2u; -N\delta, u-c+\delta+i\delta,
u+d+(N-1)\delta-i\delta, 
\nonumber\\
&\qquad \qquad\qquad\qquad\qquad
u+a+j\delta,u-b+\delta-j\delta,u+e_1,\ldots,u+e_r\Big|\dfrac{w}{z}\Big], 
\\
\Psi_{ij}&=
\dfrac{[c\pm u]_j}{[d\pm u]_j}
{}_{r+10}V_{r+9}\Big[2u; -N\delta, u-a+\delta+i\delta,
u+b+(N-1)\delta-i\delta, 
\nonumber\\
&\qquad \qquad\qquad\qquad\qquad
u\!+\!c\!+\!j\delta,u\!-\!d\!+\!\delta\!-\!j\delta,
u\!-\!e_1\!+\!\delta,\ldots,u\!-\!e_r\!+\!\delta
\Big|\dfrac{z}{w}\Big]. 
\end{align}  

\par\medskip 
Theorem \ref{thm:EHG2} can be proved by a procedure similar to the one 
we used in the previous section.  In this case we define the matrix 
$L=\big(L_{ij}\big)_{i,j=0}^N$ by
\begin{align}
L_{ij}=V^{(j)}\Big[
2u;\, -N\delta, u-c+(1+i)\delta, u+d+(N-1-i)\delta, u+c, u-d+\delta
\Big]
\end{align}
for $0\le i,j\le N$.  Then one can show
\begin{align}
(LG)_{ij}=\dfrac{[c+d+(j-i-1)\delta]_N[-(i+j)\delta]_N[2u+\delta]_N[d-c]_N}
{[u+d+j\delta]_N[u-c+(1-j)\delta]_N[u+c-i\delta]_N[-u+d-(1+i)\delta]_N}
\dfrac{[c\pm u]_j}{[d\pm u]_j}
\end{align}
by means of the Frenkel-Turaev sum, and hence 
$(LG)_{ij}=0$ for $i+j<N$.  
Then the series $\Phi_{ij}$ are obtained by computing the product 
$L\widetilde{F}$ as before.  
We remark that in this case
\begin{align}
\det M&=(-1)^{\binom{n+1}{2}}\dprod{j=0}{n}\,\dfrac{[c\pm u]_j}{[d\pm u]_j}
\nonumber\\
&\quad\cdot 
\dprod{j=0}{n}\,
\dfrac{[c+d-(N+1-2j)\delta]_N[-N\delta]_N[2u+\delta]_N[d-c]_N}
{[u\!+\!d\!+\!j\delta]_N[u\!-\!c\!+\!(1\!-\!j)\delta]_N
[u\!+\!c\!-\!(N\!-\!j)\delta]_N[-u\!+\!d\!-\!(N\!+\!1\!-\!j)\delta]_N}. 
\end{align}
The determinant of $L$ can also be computed 
in a factorized form by the elliptic version \eqref{eq:B14} of Krattenthaler's formula:  
\begin{align}
\det L&=
\dfrac{
\tprod{j=1}{N}[\delta]_j[c+d+(N-1-2j)\delta]_j[c-d-(N-1)\delta]_j[2u+j\delta]_j
}{
\tprod{i=0}{N} [u+c-i\delta]_N[u-d-(N-2-i)\delta]_N
}
\nonumber\\
&\quad\cdot
\dprod{j=0}{N}
\dfrac{[2u+2j\delta]}{[2u]}
\dfrac{
[2u]_j[-N\delta]_j[u+c]_j[u-d+\delta]_j
}{
[\delta]_j[2u+(N+1)\delta]_j[u-c+\delta]_j[u+d]_j
}.  
\end{align}
The constant in \eqref{eq:PDetEHG2} is given by $K_{m,n}(c;d)=\det M/\det L$.  

\section*{Acknowledgment}
The author would like to express his thanks to the anonymous referee 
for providing various informations on preceding works relevant to 
the subject of this paper.

\appendix
\section{Condensation of determinants} 

In this Appendix A, we give a review on the variation of Dodgson condensation 
(Sylvester identity) of determinants 
due to Gasca--L\'opez-Carmona--Ramirez \cite{GLR1982},
which we call the {\em condensation along a moving core}. 
For further generalizations of Sylvester's identity, we refer the reader 
to M\"uhlbach--Gasca \cite{MG1985}.  

\par\medskip
We first recall a standard version of the Dodgson condensation (Sylvester's identity)
for comparison. 
For a general $m\times n$ matrix 
$X=\big(x_{ij}\big)_{1\le i\le m, 1\le j\le n}$ 
(with entries in a commutative ring), 
we denote by 
$X^{i_1,\ldots,i_r}_{j_1,\ldots,j_s}=\big(x_{i_aj_b}\big)_{1\le a\le r,1\le b\le s}$ 
the $r\times s$ submatrix of $X$ with row 
indices $i_1,\ldots,i_r\in\{1,\ldots,m\}$ and 
column indices $j_1,\ldots,j_r\in\{1,\ldots,n\}$.  
When $r=s$, we denote by $\det X^{i_1,\ldots,i_r}_{j_1,\ldots,j_r}$ the corresponding 
minor determinant. 

\begin{lem}[Dodgson condensation, Sylvester's identity]\label{lem:ACond1}
Let 
$X=\big(x_{ij}\big)_{i,j=1}^n$ 
an $n\times n$ matrix and set $n=r+s$ $(r,s\ge 1)$. 
%With the notation of minor determinants 
%$\det X^{i_1,\ldots,i_l}_{j_1,\ldots,j_l}
%=\det\left(x^{i_a}_{j_b}\right)_{a,b=1}^l$, 
We define an $r\times r$ matrix $Y=\big(y_{ij}\big)_{i,j=1}^r$ by 
using the $(s+1)\times(s+1)$ minor determinants 
$y_{ij}=\det X^{i,r+1,\ldots,n}_{j,r+1,\ldots,n}$ of $X$. 
%\begin{equation}
%Y=(y^{i}_j)_{i,j=1}^{r};\qquad y^{i}_j=\det X^{i,i+1,\ldots,i+s}_{j,r+1,\ldots,n}. 
%\end{equation}
Then the determinant of $Y$ is expressed as
\begin{equation}\label{eq:Dodgson}
\det Y 
=\det{X}\ (\det X^{r+1,\ldots,n}_{r+1,\ldots,n})^{r-1};
\quad
Y=\big(y_{ij}\big)_{i,j=1}^{r},\ \ y_{ij}=\det X^{i,r+1,\ldots,n}_{j,r+1,\ldots,n}. 
\end{equation}
\end{lem}
\proof{}{
Define an $n\times n$ upper triangular matrix $Z=\big(z_{ij}\big)_{i,j=1}^n$ 
by setting 
\begin{align}
z_{ij}=
\left\{
\begin{array}{ll}
\delta_{i,j}\det X^{r+1\ldots,n}_{r+1,\ldots,n}\quad&(1\le i,j\le r)\\
(-1)^{j-r}\det X^{i,r+1,\ldots,\widehat{j},\ldots,n}_{r+1,\,\ldots,\,n}\quad
&(1\le i\le r;\ r+1\le j\le n)\\[4pt]
\delta_{ij} &(\mbox{otherwise}).  
\end{array}
\right.
\end{align}
Then for $i=1,\ldots,r$, the $(i,j)$-component of the product $ZX$ is given by 
\begin{align}
(ZX)_{ij}&=z_{ii}\,x_{ij}+\dsum{k=r+1}{}z_{ik}\,x_{kj}
\nonumber\\
&=\det X^{r+1,\ldots,n}_{r+1,\ldots,n}\,x_{ij}+\dsum{k=r+1}{n}(-1)^{k-r}
\det X^{r+1,\ldots,\widehat{k},\ldots,n}_{r+1,\ldots,n}\,x_{kj}
\nonumber\\
&=\det X^{i,r+1,\ldots,n}_{j,r+1,\ldots,n}.
\end{align}
namely,
\begin{align}
(ZX)_{ij}=\left\{
\begin{array}{ll}
y_{ij}\quad& (j=1,\ldots,r)\\
0 & (j=r+1,\ldots,n)
\end{array}
\right.
\end{align}
This means that 
\begin{align}
\arraycolsep=2pt
Z X=\left[
\begin{matrix}
Y & 0 \\[4pt]
X^{r+1,\ldots,n}_{1,\ldots,r} & 
X^{r+1,\ldots,n}_{r+1,\ldots,n}
\end{matrix}
\right]. 
\end{align}
Since $\det Z=z_{11}\cdots z_{rr}$, we obtain
\begin{align}
\det X \ (\det X^{r+1,\ldots,n}_{r+1,\ldots,n})^r
=\det Y\,\det X^{r+1,\ldots,n}_{r+1,\ldots,n}.  
\end{align}
This implies the polynomial identity
\begin{align}
\det Y = \det X \ (\det X^{r+1,\ldots,n}_{r+1,\ldots,n})^{r-1}
\end{align}
in the variables $x_{ij}$ ($1\le i,j\le n$). 
}

When $(r,s)=(n-1,1)$, \eqref{eq:Dodgson} means that
\begin{align}
\det\big(x_{ij}x_{nn}-x_{in}x_{nj}\big)_{i,j=1}^{n-1}=\det X\,x_{nn}^{n-2}.  
\end{align}
Another extreme case $(r,s)=(2,n-2)$ implies  
\begin{align}\label{eq:Jacobi}
\det X^{1,3,\ldots,n}_{1,3,\ldots,n}\det X^{2,3,\ldots,n}_{2,3,\ldots,n}
-
\det X^{1,3,\ldots,n}_{2,3,\ldots,n}\det X^{2,3,\ldots,n}_{1,3,\ldots,n}
=\det X^{1,2,\ldots,n}_{1,2,\ldots,n}\,\det X^{3,\ldots,n}_{3,\ldots,n}, 
\end{align}
which is equivalent to 
\begin{align}\label{eq:Lewis-Carroll}
\det X^{1,\ldots,n-1}_{1,\ldots,n-1}\det X^{2,\ldots,n}_{2,\ldots,n}
-
\det X^{1,\ldots,n-1}_{2,\ldots,n}\det X^{2,\ldots,n}_{1,\ldots,n-1}
=\det X^{1,2,\ldots,n}_{1,2,\ldots,n}\,\det X^{2.\ldots,n-1}_{2,\ldots,n-1}. 
\end{align}
These identities \eqref{eq:Jacobi}, \eqref{eq:Lewis-Carroll} are 
often referred to as {\em Jacobi's formula} or {\em Lewis--Carroll's formula}. 

\par\medskip
The variant of Dodgson condensation that we use in this paper 
is the following identity due to Gasca--L\'opez-Carmona--Ramirez \cite{GLR1982}. 

\begin{lem}[Condensation along a moving core] \label{lem:ACond2}
Let 
$X=\big(x_{ij}\big)_{i,j=1}^n$ an $n\times n$ matrix and set $n=r+s$ $(r,s\ge 1)$. 
%With the notation of minor determinants 
%$\det X^{i_1,\ldots,i_l}_{j_1,\ldots,j_l}
%=\det\left(x^{i_a}_{j_b}\right)_{a,b=1}^l$, 
We define an $r\times r$ matrix $Y=\big(y_{ij}\big)_{i,j=1}^r$ by 
using the $(s+1)\times(s+1)$ minor determinants 
$y_{ij}=\det X^{i,i+1,\ldots,i+s}_{j,r+1,\ldots,n}$ of $X$. 
%\begin{equation}
%Y=(y^{i}_j)_{i,j=1}^{r};\qquad y^{i}_j=\det X^{i,i+1,\ldots,i+s}_{j,r+1,\ldots,n}. 
%\end{equation}
Then the determinant of $Y$ is expressed as
\begin{equation}\label{eq:ACond2}
\det Y 
=\det{X}\ \dprod{i=1}{r-1}\,\det X^{i+1,\ldots,i+s}_{r+1,\ldots,n}; 
\quad
Y=\big(y_{ij}\big)_{i,j=1}^{r},\ \ y_{ij}=\det X^{i,i+1,\ldots,i+s}_{j,r+1,\ldots,n}. 
\end{equation}
\end{lem}
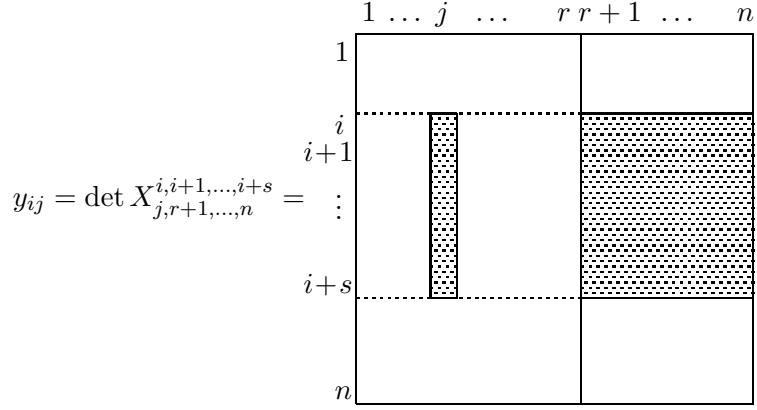
\begin{figure}[h]
$$
\begin{picture}(320,165)(-120,0)
\put(-130,85){$y_{ij}=\det X^{i,i+1,\ldots,i+s}_{j,r+1,\ldots,n}=$}
\put(0,10){\line(1,0){150}}
\put(0,150){\line(1,0){150}}
\put(0,10){\line(0,1){140}}
\put(150,10){\line(0,1){140}}
\put(85,10){\line(0,1){140}}
\put(2,155){$1$}
\put(12,155){$\ldots$}
\put(30,155){$j$}
\put(45,155){$\ldots$}
\put(76,155){$r$}
\put(84,155){$r+1$}
\put(115,155){$\ldots$}
\put(144,155){$n$}
\put(85,120){\line(1,0){65}}
\put(85,50){\line(1,0){65}}
\multiput(0,120)(3,0){30}{\line(1,0){1}}
\multiput(0,50)(3,0){30}{\line(1,0){1}}
\put(28,120){\line(1,0){10}}
\put(28,50){\line(1,0){10}}
\put(28,50){\line(0,1){70}}
\put(38,50){\line(0,1){70}}
\put(-8,140){$1$}
\put(-8,112){$i$}
\put(-20,102){$i\!+\!1$}
\put(-8,80){$\vdots$}
\put(-20,52){$i\!+\!s$}
\put(-8,12){$n$}
\multiput(28,118)(3,0){4}{\line(1,0){1}}
\multiput(85,118)(3,0){22}{\line(1,0){1}}
\multiput(29.5,116)(3,0){3}{\line(1,0){1}}
\multiput(86.5,116)(3,0){22}{\line(1,0){1}}
\multiput(28,114)(3,0){4}{\line(1,0){1}}
\multiput(85,114)(3,0){22}{\line(1,0){1}}
\multiput(29.5,112)(3,0){3}{\line(1,0){1}}
\multiput(86.5,112)(3,0){22}{\line(1,0){1}}
\multiput(28,110)(3,0){4}{\line(1,0){1}}
\multiput(85,110)(3,0){22}{\line(1,0){1}}
\multiput(29.5,108)(3,0){3}{\line(1,0){1}}
\multiput(86.5,108)(3,0){22}{\line(1,0){1}}
\multiput(28,106)(3,0){4}{\line(1,0){1}}
\multiput(85,106)(3,0){22}{\line(1,0){1}}
\multiput(29.5,104)(3,0){3}{\line(1,0){1}}
\multiput(86.5,104)(3,0){22}{\line(1,0){1}}
\multiput(28,102)(3,0){4}{\line(1,0){1}}
\multiput(85,102)(3,0){22}{\line(1,0){1}}
\multiput(29.5,100)(3,0){3}{\line(1,0){1}}
\multiput(86.5,100)(3,0){22}{\line(1,0){1}}
\multiput(28,98)(3,0){4}{\line(1,0){1}}
\multiput(85,98)(3,0){22}{\line(1,0){1}}
\multiput(29.5,96)(3,0){3}{\line(1,0){1}}
\multiput(86.5,96)(3,0){22}{\line(1,0){1}}
\multiput(28,94)(3,0){4}{\line(1,0){1}}
\multiput(85,94)(3,0){22}{\line(1,0){1}}
\multiput(29.5,92)(3,0){3}{\line(1,0){1}}
\multiput(86.5,92)(3,0){22}{\line(1,0){1}}
\multiput(28,90)(3,0){4}{\line(1,0){1}}
\multiput(85,90)(3,0){22}{\line(1,0){1}}
\multiput(29.5,88)(3,0){3}{\line(1,0){1}}
\multiput(86.5,88)(3,0){22}{\line(1,0){1}}
\multiput(28,86)(3,0){4}{\line(1,0){1}}
\multiput(85,86)(3,0){22}{\line(1,0){1}}
\multiput(29.5,84)(3,0){3}{\line(1,0){1}}
\multiput(86.5,84)(3,0){22}{\line(1,0){1}}
\multiput(28,82)(3,0){4}{\line(1,0){1}}
\multiput(85,82)(3,0){22}{\line(1,0){1}}
\multiput(29.5,80)(3,0){3}{\line(1,0){1}}
\multiput(86.5,80)(3,0){22}{\line(1,0){1}}
\multiput(28,78)(3,0){4}{\line(1,0){1}}
\multiput(85,78)(3,0){22}{\line(1,0){1}}
\multiput(29.5,76)(3,0){3}{\line(1,0){1}}
\multiput(86.5,76)(3,0){22}{\line(1,0){1}}
\multiput(28,74)(3,0){4}{\line(1,0){1}}
\multiput(85,74)(3,0){22}{\line(1,0){1}}
\multiput(29.5,72)(3,0){3}{\line(1,0){1}}
\multiput(86.5,72)(3,0){22}{\line(1,0){1}}
\multiput(28,70)(3,0){4}{\line(1,0){1}}
\multiput(85,70)(3,0){22}{\line(1,0){1}}
\multiput(29.5,68)(3,0){3}{\line(1,0){1}}
\multiput(86.5,68)(3,0){22}{\line(1,0){1}}
\multiput(28,66)(3,0){4}{\line(1,0){1}}
\multiput(85,66)(3,0){22}{\line(1,0){1}}
\multiput(29.5,64)(3,0){3}{\line(1,0){1}}
\multiput(86.5,64)(3,0){22}{\line(1,0){1}}
\multiput(28,62)(3,0){4}{\line(1,0){1}}
\multiput(85,62)(3,0){22}{\line(1,0){1}}
\multiput(29.5,60)(3,0){3}{\line(1,0){1}}
\multiput(86.5,60)(3,0){22}{\line(1,0){1}}
\multiput(28,58)(3,0){4}{\line(1,0){1}}
\multiput(85,58)(3,0){22}{\line(1,0){1}}
\multiput(29.5,56)(3,0){3}{\line(1,0){1}}
\multiput(86.5,56)(3,0){22}{\line(1,0){1}}
\multiput(28,54)(3,0){4}{\line(1,0){1}}
\multiput(85,54)(3,0){22}{\line(1,0){1}}
\multiput(29.5,52)(3,0){3}{\line(1,0){1}}
\multiput(86.5,52)(3,0){22}{\line(1,0){1}}
\multiput(28,50)(3,0){4}{\line(1,0){1}}
\multiput(85,50)(3,0){22}{\line(1,0){1}}
\end{picture}
$$
\caption{Condensation along a moving core}
\end{figure}
% We leave the proof of this lemma to the reader as an exercise of linear algebra.  
\proof{}{
We define an $n\times n$ upper triangular 
matrix $Z=\big(z_{ij}\big)_{i,j=1}^{n}$ as 
follows by using $s\times s$ minor determinants of $X$:  
\begin{equation}
z_{ij}=
\left\{
\begin{array}{ll}
(-1)^{j-i}\det X^{i,\ldots,\widehat{j},\ldots,i+s}_{r+1,\,\ldots,\,n}\qquad
&(1\le i\le r; i\le j\le i+s)\\[6pt]
\delta_{ij}&(\mbox{otherwise}).
\end{array}
\right.
\end{equation}
Then for $i=1,\ldots,r$, we have
\begin{align}
\big(ZX\big)_{ij}=\dsum{k=i}{i+s}z_{ik}\,x_{kj}
=\dsum{k=i}{i+s}(-1)^{k-i}\det X^{i,i+1,\ldots,\widehat{k},\ldots,i+s}_{r+1,\,\ldots,\,n}
\,x_{kj}
=\det X^{i,\ldots,i+s}_{j,r+1,\ldots,n},
\end{align}
namely
\begin{align}
\big(ZX\big)_{ij}=\left\{\begin{array}{ll}
y_{ij} \quad&(1\le j\le r)
\\[6pt]
0 & (r+1\le j\le n).
\end{array}\right.
\end{align}
This means that 
\begin{align}
\arraycolsep=2pt
Z X=\left[
\begin{matrix}
Y & 0 \\[4pt]
X^{r+1,\ldots,n}_{1,\ldots,r} & 
X^{r+1,\ldots,n}_{r+1,\ldots,n}
\end{matrix}
\right]. 
\end{align}
Since $\det Z=z_{11}\cdots z_{rr}$, we obtain
\begin{align}
\det X \ \dprod{i=1}{r}\,\det X^{i,\ldots,i+s}_{r+1,\ldots,n}
=\det Y\,\det X^{r+1,\ldots,n}_{r+1,\ldots,n}.  
\end{align}
This implies the polynomial identity
\begin{align}
\det Y = \det X \ \dprod{i=1}{r-1}\,\det X^{i,\ldots,i+s}_{r+1,\ldots,n}
\end{align}
in the variables $x_{ij}$ ($1\le i,j\le n$). 
}

We remark that, 
if we renormalize the matrix $Y$ by setting 
\begin{align}\label{eq:A18}
\widetilde{Y}=\big(\widetilde{y}_{ij}\big)_{i,j=1}^{r},\quad
\widetilde{y}_{ij}=
\det X^{i,i+1,\ldots,i+s}_{j,r+1,\ldots,n}\,
\big(\det X^{i+1,\ldots,i+s}_{r+1,\ldots,n}\big)^{-1}
\quad(i,j=1,\ldots,r), 
\end{align}
then equality \eqref{eq:ACond2} is rewritten equivalently as
\begin{align}\label{eq:A19}
%\det \widetilde{Y}=\det X\, \big(\det X^{r+1,\ldots,n}_{r+1,\ldots,n}\big)^{-1},
%\quad\mbox{namely,}\quad
\det X=
\det \widetilde{Y} \det X^{r+1,\ldots,n}_{r+1,\ldots,n}.  
\end{align}

\section{Variations of Krattenthaler's determinant formula}

In this Appendix B, we recall Krattenthaler's determinant formula \cite{K1999} 
and its elliptic extension due to Warnaar \cite{W2002}.  
Although these formulas can be proved in various ways, 
we remark here that they are consequences of 
Lemma \ref{lem:abstractK} below, which can be regarded as an 
abstract form of Krattenthaler's determinant formula
(for recent works on the evaluation of determinants involving 
shifted factorials, see Normand \cite{N2004}).

\par\medskip
We first recall a typical form of Krattenthaler's determinant formula \cite{K1999}. 
\begin{lem}\label{lem:Krattenthaler}
For any set of variables $x_i$ $(0\le i\le m)$ and parameters
$\alpha_k,\beta_k,\gamma_k,\delta_k$ $(0\le k<m)$, 
one has 
\begin{align}\label{eq:Krattenthaler}
\det\left(\dprod{0\le k<j}{}\dfrac{\alpha_k x_i+\beta_k}
{\gamma_k x_i+\delta_k}\right)_{i,j=0}^{m}
=
\dfrac{
\dprod{0\le i<j\le m}{}(x_j-x_i)
\dprod{0\le k\le l<m}{}(\alpha_k\delta_l-\beta_k\gamma_l)
}
{\dprod{0\le i\le m}{}\,\dprod{0\le k<m}{} (\gamma_kx_i+\delta_k)}. 
\end{align}
\end{lem}
By specializing the parameters $\alpha_k,\beta_k,\gamma_k,\delta_k$, 
we obtain various determinant formulas.  We quote below some of them. 
\par\smallskip
\noindent
(a)\quad Case where $\alpha_k=\gamma_k=1$ and $\beta_k=a_k$, $\delta_k=b_k$: 
\begin{align}
\det\left(\dprod{0\le k<j}{}\dfrac{x_i+a_k}
{x_i+b_k}\right)_{i,j=0}^{m}
=
\dfrac{
\dprod{0\le i<j\le m}{}(x_j-x_i)
\dprod{0\le k\le l<m}{}(b_l-a_k)
}
{\dprod{0\le i\le m}{}\,\dprod{0\le k<m}{} (x_i+b_k)}. 
\end{align}
In particular, by setting $a_k=a+k$, $b_k=b+k$ one obtains
\begin{align}\label{eq:B3}
\det\left(
\dfrac{(a+x_i)_j}{(b+x_i)_j}
\right)_{i,j=0}^{m}
=
\dfrac{
\dprod{0\le i<j\le m}{}(x_j-x_i)\,
\dprod{k=1}{m}\,(b-a)_k
}
{\dprod{0\le i\le m}{}\,(b+x_i)_m}. 
\end{align}
where $(a)_k=a(a+1)\cdots(a+k-1)$.  
\newline
(b)\quad Case where $\beta_k=\delta_k=1$, $\alpha_k=-a_k$, $\gamma_k=-b_k$: 
\begin{align}
\det\left(\dprod{0\le k<j}{}\dfrac{1-a_k x_i}
{1-b_kx_i}\right)_{i,j=0}^{m}
=
\dfrac{
\dprod{0\le i<j\le m}{}(x_j-x_i)
\dprod{0\le k\le l<m}{}(b_l-a_k)
}
{\dprod{0\le i\le m}{}\,\dprod{0\le k<m}{} (1-b_kx_i)}. 
\end{align}
By setting $a_k=p^ka$ and $b_k=q^kb$, one has
\begin{align}
\det\left(\dprod{0\le k<j}{}\dfrac{(ax_i;p)_j}
{(bx_i;q)_j}\right)_{i,j=0}^{m}
=
\dfrac{
\dprod{0\le i<j\le m}{}(x_j-x_i)
\dprod{0\le k\le l<m}{}(q^lb-p^ka)
}
{\dprod{0\le i\le m}{}\, (bx_i;q)_m}. 
\end{align}
where $(a;p)_k=(1-a)(1-pa)\cdots(1-p^{k-1}a)$.  
In particular, 
\begin{align}
\det\left(\dprod{0\le k<j}{}\dfrac{(ax_i;q)_j}
{(bx_i;q)_j}\right)_{i,j=0}^{m}
=
\dfrac{a^{\binom{m+1}{2}}q^{\binom{m+1}{3}}
\dprod{0\le i<j\le m}{}(x_i-x_j)
\dprod{k=1}{m}(b/a;q)_k
}
{\dprod{0\le i\le m}{}\,(bx_i;q)_m}. 
\end{align}
(c)\quad Noting that 
\begin{align}
(az;q)_n(ac/z;q)_n=\dprod{0\le k<n}{}(1+a^2cq^{2k}-aq^k(z+c/z))
\end{align}
set
\begin{align}
x_i=z_i+c/z_i;\quad \alpha_k=-ap^{k},\ \ \beta_k=1+a^2c p^{2k},\ \ 
\gamma_k=-b q^{k},\ \ \delta_k=1+b^2 cq^{2k}. 
\end{align}
Then we have 
\begin{align}
&\det\left(
\dfrac{(az_i;p)_j(ac/z_i;p)_j}{(bz_i;q)_j(bc/z_i;q)_j}
\right)_{i,j=0}^{m}
\nonumber\\
&
=
\dfrac{
\dprod{0\le i<j\le m}{}(z_j-z_i)(1-c/z_iz_j)
\dprod{0\le k\le l<m}{}(bq^l-ap^k)(1-p^kq^labc)
}
{
\dprod{0\le i\le m}{}\,(bz_i;q)_m(bc/z_i;q)_m
}. 
\end{align}
In particular,
\begin{align}
&\det\left(
\dfrac{(az_i;q)_j(ac/z_i;q)_j}{(bz_i;q)_j(bc/z_i;q)_j}
\right)_{i,j=0}^{m}
\nonumber\\
&
=
\dfrac{
a^{\binom{m+1}{2}}q^{\binom{m+1}{3}}
\dprod{0\le i<j\le m}{}(z_i-z_j)(1-c/z_iz_j)\ 
\dprod{k=1}{m}(b/a;q)_k(q^{2(m-1-k)}abc;q)_k
}
{
\dprod{0\le i\le m}{}\,(bz_i;q)_m(bc/z_i;q)_m
}. 
\end{align}

\par\medskip
Let $[x]$ a nonzero entire function in $x\in\mathbb{C}$ and 
suppose that $[x]$ satisfies the so-called Riemann relation:  
For any $x,\alpha,\beta,\gamma\in\mathbb{C}$, 
\begin{equation}
[x\pm\alpha][\beta\pm\gamma]
+
[x\pm\beta][\gamma\pm\alpha]
+
[x\pm\gamma][\alpha\pm\beta]
=0,
\end{equation}
where $[x\pm \alpha]=[x+\alpha][x-\alpha]$.  
This functional equation is equivalent to 
\begin{align}
[x\pm u][y\pm v]-[x\pm v][y\pm u]=[x\pm y][u\pm v].  
\end{align}
The following lemma is the elliptic extension of Lemma \ref{lem:Krattenthaler}
due to Warnaar \cite{W2002}.   
\begin{lem}\label{lem:Warnaar}
For any set of variables $x_i$ $(0\le i\le m)$ and parameters
$a_k, b_k$ $(0\le k<m)$, 
one has 
\begin{align}\label{eq:Warnaar}
\det\left(\dprod{0\le k<j}{}
\dfrac{[a_k\pm x_i]}{[b_k\pm x_i]}
\right)_{i,j=0}^{m}
=
\dfrac{
\dprod{0\le i<j\le m}{}[x_j\pm x_i]
\dprod{0\le k\le l<m}{}[a_k\pm b_l]
}
{\dprod{0\le i\le m}{}\,\dprod{0\le k<m}{} [b_k\pm x_i]}. 
\end{align}
\end{lem}
As a special case where 
$a_k=a+k\delta$, $b_k=b+k\delta$ ($0\le k<m$), we obtain 
\begin{align}\label{eq:B14}
\det\left(
\dfrac{[a\pm x_i]_j}{[b\pm x_i]_j}
\right)_{i,j=0}^{m}
=
\dfrac{
\dprod{0\le i<j\le m}{}[x_i\pm x_j]\,
\dprod{k=1}{m}[b-a]_k[a+b+(k-1)\delta]_k
}
{\dprod{0\le i\le m}{}\,[b\pm x_i]_m}. 
\end{align}
where $[a]_k=[a][a+\delta]\cdots[a+(k-1)\delta]$ 
and $[a\pm b]_k=[a+b]_k[a-b]_k$.  

\par\medskip 
Lemma \ref{lem:Krattenthaler} and Lemma \ref{lem:Warnaar} can be 
proved as consequences of the following abstract form of 
Krattenthaler's determinant formula.

Let $a_{ik}, b_{ik}$ ($0\le i\le N$; $0\le k<N$) be elements 
a field $\mathbb{K}$ with $b_{ik}\ne 0$ for all $i,k$, and consider 
the matrix 
\begin{align}
X_m=\Bigg(\dprod{0\le k<j}{}\dfrac{a_{ik}}{b_{ik}}\Bigg)_{i,j=0}^m.
\end{align}
for $m=0,1,\ldots,N$. 
Suppose that there exist elements $p_{ij}$ ($0\le i,j\le N$), 
$q_{kl}$ ($0\le k\le l<N$) of $\mathbb{K}$ such that
\begin{align}\label{eq:factorization}
a_{ik}b_{jl}-a_{jk}b_{il}=p_{ij}q_{kl},\quad p_{ji}=-p_{ij}
\end{align}
for all $i,j\in\{0,1\ldots,N\}$ and $k,l\in\{0,1,\ldots,N-1\}$. 
\begin{lem}\label{lem:abstractK}
Under the assumption \eqref{eq:factorization}, the determinant $\det X_m$ 
is factorized as 
\begin{align}\label{eq:abstractK}
\det X_m=\det \Bigg(\dprod{0\le k<j}{}\dfrac{a_{ik}}{b_{ik}}\Bigg)_{i,j=0}^m
=\dfrac{\dprod{0\le i<j\le m}{}p_{ji}\,\dprod{0\le k\le l<m}{}q_{kl}}
{\dprod{i=0}{m}\dprod{k=0}{m-1} b_{i,k}}.  
\end{align}
for $m=0,1,\ldots,N$.
\end{lem}
Set $\tau_m=\det X_m$ for $m=0,1,2,\ldots,N$, so that 
\begin{align}
\tau_0=1,\quad  \tau_1=\det\left[\begin{matrix}
1 & \dfrac{a_{00}}{b_{00}}\\[10pt]
1 & \dfrac{a_{10}}{b_{10}}
\end{matrix}\right]
=\dfrac{a_{10}}{b_{10}}-\dfrac{a_{00}}{b_{00}}. 
\end{align}
The first nontrivial case is 
guaranteed by the assumption \eqref{eq:factorization}: 
\begin{align}
\tau_1
%=\dfrac{a_{10}}{b_{10}}-\dfrac{a_{00}}{b_{00}}
=\dfrac{a_{10}b_{00}-a_{00}b_{10}}{b_{00}b_{10}}
=\dfrac{p_{10}\,q_{00}}{b_{00}\,b_{10}}.  
\end{align}
Then \eqref{eq:abstractK} can 
be proved by means of the Lewis--Carroll formula.  
In fact from \eqref{eq:Lewis-Carroll}, we obtain the bilinear identities 
\begin{align}
\dfrac{a_{m+1,1}}{b_{m+1,1}}\,\tau_m\,T_CT_R(\tau_m)
-\dfrac{a_{11}}{b_{11}}\,T_C(\tau_m)\,T_R(\tau_m)
=\tau_{m+1}\,T_CT_R(\tau_{m-1}),
\end{align}
for $\tau_m$, 
where $T_R$ and $T_C$ stand for the symbolic shift operator for the 
row indices and the column indices:
\begin{align}\label{eq:bilintau}
&T_{R}(a_{ij})=a_{i+1,j},\ \ T_{R}(b_{ij})=b_{i+1,j},\ \ 
T_{R}(p_{ij})=p_{i+1,j+1},\ \ T_{R}(q_{kl})=q_{kl},
\nonumber\\
&T_{C}(a_{ij})=a_{i,j+1},\ \ T_{C}(b_{ij})=b_{i,j+1},\ \ 
T_{C}(p_{ij})=p_{ij},\ \ T_{C}(q_{kl})=q_{k+1,l+1}.  
\end{align}
Thanks to the bilinear identities, one can prove 
\begin{align}
\tau_m=\dfrac{\dprod{0\le i<j\le m}{}p_{ji}\,\dprod{0\le k\le l<m}{}q_{kl}}
{\dprod{i=0}{m}\dprod{k=0}{m-1} b_{i,k}}
\qquad(m=0,1,\ldots,N)
\end{align}
by the induction on $m$. 
\par\medskip
We remark that Lemma \ref{lem:Krattenthaler} is the case 
where $a_{ik}=\alpha_kx_i+\beta_k$, $b_{ik}=\gamma_k x_i+\delta_k$. 
Since 
\begin{align}\label{eq:rational}
(\alpha_kx_i+\beta_k)(\gamma_l x_j+\delta_l)
-(\alpha_kx_j+\beta_k)(\gamma_l x_i+\delta_l)
=(x_i-x_j)(\alpha_k\delta_l-\beta_k\gamma_l)
\end{align}
the factorization condition \eqref{eq:factorization} is verified 
with $p_{ij}=x_i-x_j$ and $q_{kl}=\alpha_k\delta_l-\beta_k\gamma_l$. 
Lemma \ref{lem:Warnaar} is the case where 
$a_{ik}=[a_k\pm x_i]$, $b_{ik}=[b_k\pm x_i]$. 
Since
\begin{align}
[a_k\pm x_i][b_l\pm x_j]-
[a_k\pm x_j][b_l\pm x_i]=[a_k\pm b_l][x_i\pm x_j]
\end{align}
the condition \eqref{eq:factorization} is satisfied with 
$p_{ij}=[x_i\pm x_j]$ and $q_{kl}=[a_k\pm b_l]$. 
One can prove in fact that {\em generic} solutions to the system of 
equations \eqref{eq:factorization} reduce to the case of \eqref{eq:rational}.  
It would be worthwhile, however, to recognize the role of bilinear equations 
\eqref{eq:bilintau} which lead to the factorization of determinants. 

%%%%%%%%%%%%%%%%%%%%%%%%%%%%%%%%%%%%

%%%%%%%%%%%%%%%%%%%%%%%%%%%%%%%%%%%%

\end{document}